\RequirePackage{fix-cm}
\documentclass[smallextended,numbook,envcountsame]{svjour3M}       
\smartqed  
\usepackage{graphicx}
\usepackage{pdfsync}

\usepackage{epsfig}
\usepackage{mathptmx}
\usepackage{amsmath,amsfonts,amssymb}
\usepackage{mathtools}
\usepackage{mathrsfs}
\usepackage[all]{xy}
\usepackage{graphicx}
\usepackage{latexsym}
\usepackage{verbatim}

\usepackage{tikz}
\usetikzlibrary{calc}
\usetikzlibrary{arrows.meta}

\numberwithin{equation}{section}

\renewcommand{\Re}{\mathop{\mathrm{Re}}}
\renewcommand{\Im}{\mathop{\mathrm{Im}}}
\newcommand{\id}{\mathsf{id}}
\newcommand{\Hol}{\mathsf{Hol}}

\renewcommand{\le}{\leqslant}
\renewcommand{\ge}{\geqslant}
\renewcommand{\leq}{\leqslant}
\renewcommand{\geq}{\geqslant}

\newcommand{\interior}{\mathop{\mathsf{int}}}

\newcommand{\Arg}{\mathop{\mathsf{Arg}}}

\newcommand{\Real}{\mathbb{R}}
\newcommand{\Natural}{\mathbb{N}}
\newcommand{\N}{\mathbb{N}}
\newcommand{\Complex}{\mathbb{C}}

\newcommand{\UD}{\mathbb{D}}

\newcommand{\D}{\mathbb D}

\newcommand{\C}{\mathbb C}
\newcommand{\R}{\mathbb R}
\newcommand{\Aut}{\mathsf{Aut}(\mathbb D)}
\newcommand{\de}{\partial}
\newcommand{\dist}{\,\mathop{\mathsf{dist}\!}}
\newcommand{\const}{\mathrm{const}}

\newcommand{\Ker}[2]{\mathcal{K}\big(#1,#2\big)}

\newcommand{\vstrip}{V}

\long\def\REM#1{\relax}

\begin{document}

\title{Angular extents and trajectory slopes in the theory of holomorphic semigroups in the unit disk \thanks{$^\dag$ Partially supported by the \textit{Ministerio
de Econom\'{\i}a y Competitividad} and the European Union (FEDER) PGC2018-094215-B-100 and by \textit{La Junta de Andaluc\'{\i}a}, FQM-133}
}

\author{Manuel D. Contreras\,$^\dag$         \and
        Santiago D\'{\i}az-Madrigal\,$^\dag$ \and
        Pavel Gumenyuk
}

\authorrunning{M. D. Contreras, S. D\'{\i}az-Madrigal, P. Gumenyuk}

\institute{M. D. Contreras and S. D\'{\i}az-Madrigal \at
              Camino de los Descubrimientos, s/n\\
	Departamento de Matem\'{a}tica Aplicada II and IMUS\\
	Universidad de Sevilla\\
	Sevilla 41092, SPAIN\\
              \email{contreras@us.es}\\ \email{madrigal@us.es}
\and
           P. Gumenyuk \at
              Department of
Mathematics\\ Milano Politecnico, via E. Bonardi 9\\ Milan 20133, ITALY\\
\email{pavel.gumenyuk@polimi.it}
}

\date{\today}

\maketitle

\begin{abstract}
  We study relationships between the asymptotic behaviour of a non-elliptic semigroup of holomorphic self-maps of the unit disk and the geometry of its planar domain (the image of the Koenigs function). We establish a sufficient condition for the trajectories of the semigroup to converge to its Denjoy\,--\,Wolff point with a definite slope. We obtain as a corollary two previously known sufficient conditions.
\end{abstract}

\tableofcontents

\section{Introduction}

A \textsl{(one-parameter) semigroup} $(\varphi_t)_{t\ge0}$ of holomorphic self-maps of~$\D$~--- for short,  a \textsl{semigroup in $\D$}~--- is a continuous homomorphism $t\mapsto \varphi_t$ from the
additive semigroup $(\R_{\ge0}, +)$ of non-negative real numbers to the
semigroup $\big({\sf Hol}(\D,\D),\circ\big)$ of  holomorphic self-maps
of $\D$ with respect to composition, endowed with the
topology of uniform convergence on compacta. If  ${\varphi_{t_0}}$ is an automorphism of $\D$ for some~$t_0>0$, then ${\varphi_t}$ is an automorphism for all~$t\ge0$ 
and in such a case we will say that $(\varphi_t)$ is a {\sl group}, because indeed it can be extended to a group homomorphism $\R\ni t\mapsto\varphi_t\in\Hol(\D,\D)$.

The theory of semigroups in $\D$ has a long history dating back to the early nineteen century. Moreover, nowadays, it is a flourishing branch of Analysis with strong connections with Dynamical Systems and with many applications in other areas (see \cite{BCD-Book} and the bibliography therein). Indeed, this paper is about a basic dynamical problem for semigroups in $\D$. We refer the reader to \cite{Abate}, \cite{BCD-Book}, or \cite{Shobook01} for the results  cited below without proof.

It is known that $\varphi_{t_0}$ has a fixed point in $\D$ for some $t_0>0$
if and only if there exists $\tau \in \D$ such that $\varphi_t(\tau)=\tau$  for all $t\geq 0$.
In such a case, the semigroup is called {\sl elliptic} and there exists $\lambda\in \C$
with $\Re \lambda\geq 0$ such that $\varphi_t'(\tau)=e^{-\lambda t}$ for all $t\geq 0$. The elliptic semigroup~$(\varphi_t)$ is  a group if and only if $\Re \lambda=0$. Moreover,  the above point $\tau$ is unique unless $\varphi_t=\id_\UD$ for all~${t\ge0}$, and it is called the Denjoy\,--\,Wolff point (\textsl{DW-point} in what follows) of the semigroup.

If the semigroup $(\varphi_t)$ is not elliptic, then there exists $\tau\in \partial\D$ which is the Denjoy\,--\,Wolff point of $\varphi_t$ for all $t>0$, i.e. $\varphi_t(\tau)=\tau$ and $\varphi'_t(\tau)\le1$ in the sense of angular limits. As before $\tau$ is also called  the Denjoy\,--\,Wolff point (\textsl{DW-point} in what follows) of the semigroup. In this case, there exists $\lambda\geq 0$ such that
$\varphi'_t(\tau)=e^{-\lambda t}$ for all ${t\geq 0}$, where $\varphi'_t(\tau)$ stands for the angular derivative of $\varphi_t$ at~$\tau$. A non-elliptic semigroup is said to be {\sl hyperbolic} or {\sl parabolic} depending on whether ${\lambda>0}$  or ${\lambda=0}$, respectively. Parabolic semigroups can be divided in two sub-types:
a parabolic semigroup is {\sl of positive hyperbolic step}
if $\lim_{t\to+\infty}k_\D(\varphi_{t+1}(0), \varphi_t(0))>0$, where $k_\D(\cdot,\cdot)$ denotes the hyperbolic distance in $\D$.
Otherwise, $(\varphi_t)$ is said to be  {\sl of zero hyperbolic step}.

A fundamental result for semigroups in $\D$ is the so called Continuous Denjoy\,--\,Wolff theorem which says that if $(\varphi_{t})$ is non-elliptic or elliptic but different from a group, then for any ${z\in \D}$, $\varphi_{t}(z)\to\tau$ as ${t\to+\infty}$, where $\tau$ is the Denjoy\,--\,Wolff of the semigroup. Those functions $t\mapsto \varphi_t(z)$ can be properly named orbits (or trajectories) in the usual dynamical sense thanks to  Berkson and Porta's celebrated theorem \cite[Theorem~(1.1)]{BerPor78} which asserts that
$t\mapsto \varphi_t(z)$ is real-analytic and there exists a unique
holomorphic vector field $G:\D\to \C$ such that
\[
\frac{\de \varphi_t(z)}{\de
	t}=G(\varphi_t(z)),\qquad\text{for all $z\in\D$
	and all~$t\ge0$}.
\]
This vector field $G$ is called the  {\sl vector field or infinitesimal generator} of $(\varphi_t)$.

In this paper we are interested in considering the so-called “slope problem” of the orbits of a non-elliptic semigroup in $\D$ when arriving to its Denjoy\,--\,Wolff point.

\begin{definition} Let $(\varphi_t)$ be a non-elliptic semigroup in $\D$ with Denjoy\,--\,Wolff point $\tau\in\partial\D$. The \textsl{(arrival) slope set} $\mathrm{Slope}[t\mapsto \varphi_t(z),\tau]$ of the semigroup $(\varphi_t)$ at $\tau$ with the initial point $z\in\D$ is the cluster set of the function
	$$
	[0,+\infty)\ni t\mapsto\Arg\big(1-\overline{\tau}\varphi_t(z)\big)\in (-\pi/2,\pi/2)
	$$
	as $t\to+\infty$.
	In other words, $\theta\in\left[ -\frac{\pi }{2},\frac{\pi }{2}\right]$ belongs to the set $\mathrm{Slope}[t\mapsto \varphi_t(z),\tau]$ if there exists  a sequence $(t_n)\subset[0,+\infty)$ tending to $+\infty$ such that $\Arg\big(1-\overline{\tau}\varphi_{t_n}(z)\big)\to\theta$ as ${n\to+\infty}$.
\end{definition}

\begin{remark}\label{RM_slopes} $\mathrm{Slope}[t\mapsto \varphi_t(z),\tau]$ is either a point or a closed subinterval of $\left[ -\frac{\pi }{2},\frac{\pi }{2}\right]$.
\end{remark}

For hyperbolic semigroups and parabolic semigroups of positive hyperbolic step, the arrival slope set is always a singleton (see \cite[Sect.\,17.4 and 17.5]{BCD-Book} for further information).

In contrast, for parabolic semigroups of zero hyperbolic step, the arrival slope set  does not have to reduce to a unique point (see \cite{CDG-Slope},\cite{Bet16a},\cite{Kelgiannis}). However, according to the following result by the first two authors, it does not depend on the initial point.
\begin{theopargself}
	\begin{theorem}[\protect{\cite[Theorem 2.9\,(1)]{ConDia05a}}] \label{indepen} Let $(\varphi_t)$ be a parabolic semigroup  of zero hyperbolic step with DW-point $\tau\in\partial\D$. Then, for any ${z_1,z_2\in\D}$,
		\[
		\mathrm{Slope}[t\mapsto \varphi_t(z_1),\tau]=\mathrm{Slope}[t\mapsto \varphi_t(z_2),\tau].
		\]
	\end{theorem}
\end{theopargself}

An (important) and open problem in the theory of semigroups in $\D$  has been (indeed, still is) how to detect whether the arrival slope set of a parabolic semigroup of zero hyperbolic step is a singleton or a specific kind of closed subinterval of $\left[ -\frac{\pi }{2},\frac{\pi }{2}\right]$. Here the word ``detect'' almost always means  finding sufficient and/or necessary conditions of geometric nature. This is directly related to the second key notion (with the first one being the vector field) associated with each semigroup, namely, to its holomorphic model and its Koenigs function (see \cite{BerPor78}, \cite{Cow81}, \cite{Sis}, \cite{AroBra16}, \cite[Sect.\,9]{BCD-Book}).

\begin{definition}
	Let $(\varphi_t)$ be a semigroup in $\D$. A {\sl holomorphic model} for $(\varphi_t)$ is a triple $(U, h, \Phi_t)$, where $U$ is a domain in $\C$, $(\Phi_t)$ is a group of holomorphic automorphisms of $U$, and $h: \D \to h(\D)\subset U$ is a univalent holomorphic map (called a {\sl Koenigs function} of the semigroup) satisfying the functional equation
	\begin{equation}\label{EQ_Abel0}
		h\circ \varphi_t= \Phi_t\circ h\quad\text{for all $t\ge0$}
	\end{equation}
	and the following absorbing property
	\begin{equation}\label{absorbing}
		\bigcup_{t\geq 0} \Phi_t^{-1}(h(\D))=U.
	\end{equation}
	The set $h(\D)$ is called an associated {\sl planar domain} of the semigroup.
\end{definition}

Every semigroup in $\D$ admits a holomorphic
model unique up to ``holomorphic equivalence'' (i.e., isomorphism of models). In particular, see e.g. \cite[Theorem~9.3.5 on p.\,245]{BCD-Book}, a semigroup in~$\UD$ is non-elliptic if and only if one of its (mutually equivalent) holomorphic models is of the form ${(U,h,z\mapsto z+it)}$. For such a holomorphic model, the functional equation~\eqref{EQ_Abel0} becomes Abel's classical equation
\begin{equation}\label{EQ_Abel}
h\big(\varphi_t(z)\big)=h(z)+it,\quad\text{for all~$~z\in\UD$, $t\ge0$.}
\end{equation}

Following the convention generally accepted in the literature, we will assume that all the considered holomorphic models for the non-elliptic semigroups are of the above canonical form. Then the Koenigs function becomes essentially unique:
if $h_1,h_2$ are two Koenigs functions of the same non-elliptic semigroup, then there exists a constant $c\in\C$ such that ${h_1=h_2+c}$.

Thanks to~\eqref{EQ_Abel}, planar domains of non-elliptic semigroups are complex domains of a very particular type: the so-called starlike-at-infinity domains.
\begin{definition} A domain $\Omega\subset\C$ is said to be \textsl{starlike at infinity} (in the direction of the imaginary axis) if $\Omega+it\subset \Omega$ for any~${t\geq 0}$.
\end{definition}

\begin{remark}\label{RM_Koenings-image} Any domain~$\Omega\neq\Complex$ starlike at infinity is conformally equivalent to~$\UD$ and if $h$ is a conformal mapping of~$\UD$ onto such a domain~$\Omega$, then the formula ${\varphi_t:=h^{-1}\circ(h+it)}$ for ${t\ge0}$ defines a non-elliptic semigroup in~$\UD$, whose Koenigs function is~$h$.
\end{remark}

In this context, our  problem mentioned above can be rewritten as follows: to find geometrical properties of the planar domain of a parabolic semigroup of zero hyperbolic step which imply (or characterize) whether the corresponding arrival slope set is a singleton.

As far as we know, apart from examples and some folklore results concerning strong symmetry of the planar domain, the unique three papers dealing with the above question are \cite{Bet16a}, \cite{BCDG19} and \cite{BCDGZ}.
In \cite{Bet16a}, it is shown that whenever the boundary of the planar domain is included in a vertical  half-strip, the arrival slope set is equal to $\{0\}$. Likewise, in \cite{BCDG19}, it is shown that if the boundary of the planar domain is included in a horizontal strip, the arrival slope set is also  equal to $\{0\}$.
In \cite{BCDGZ}, the authors introduces some ``boundary distance'' functions, which measure the distance of a vertical straight line to the boundary of the planar domain, and use them to characterize geometrically when the arrival slope set coincides with the singleton $\{\pi/2 \}$ or $\{-\pi/2 \}$. Moreover, they also show how these functions detect whether the convergence of the trajectories is non-tangential, i.e. whether the arrival slope set is a compact subset of ${\left( -\frac{\pi }{2},\frac{\pi }{2}\right)}$.

We would like to mention that there are also results treating the above problem in a non-geometrical way, i.e. without using planar domains. For instance, in \cite{EKRS} (see also \cite[Proposition 7.5.5]{BCD-Book}), it is proved that the arrival slope sets of a parabolic semigroup of zero hyperbolic step is a singleton whenever its vector field has enough analytic regularity (in the angular sense) at its Denjoy\,--\,Wolff point.

In this paper, we introduce some new ``angular extent'' functions of a strongly geometrical meaning, which measure the angular displacement of the boundary of the planar domain with respect to a fixed vertical straight line (see Definition \ref{angular}). Using these functions, we establish sufficient conditions for the arrival slope set of a non-elliptic semigroup to be a singleton (see Theorem~\ref{main} and Proposition~\ref{PR_second}). We also analyze the relationship between these functions and the non-tangential convergence of the orbits of the semigroup (see Proposition~\ref{prop-pendientes}). Moreover, as a corollary, we recover results from \cite{Bet16a} and \cite{BCDG19}.

The plan of the paper is as follows. In Section~2, we develop some new results about Carathéodory kernel convergence which can be of interest on their own and will be fundamental for the results of Section~5.
In Section~3, we introduce and study those angular extent functions mentioned above. Section~4 is a brief review of the boundary distance functions introduced in \cite{BCDGZ}.
We also study here their relationships with the angular extent functions from Section~3. In Section~5, we present our main results. Finally, in Section~6, we show a few examples dealing with some particularities of the angular extent functions, which, in particular,  underline important differences between them and the (apparently quite similar) boundary distance functions.

\section{Kernel convergence}
Recall the classical notion of kernel convergence of domains; for more details see e.g. \cite[\S{}II.5]{Goluzin} or \cite[\S1.4]{Pombook75}.
Let $(\Omega_n)$ be a sequence of domains in~$\C$. Fix a point $\omega\in \C$. Suppose that $\omega\in \Omega_n$ for all~$n\in\N$ large enough. Denote by $G$  the (possibly empty) set of all points $z\in \C$ possessing the following property: there exists an open connected set $\Delta\subset \C$ containing the points~$z$ and~$\omega$ and contained in~$\Omega_n$ for all sufficiently large~$n\in \N$.

The {\sl kernel} $\Ker{(\Omega_n)}{\omega}$ of the sequence~$(\Omega_n)$ with respect to the point~$\omega$ is the union ${G\cup\{\omega\}}$. The following dichotomy holds: either $G=\emptyset$ and hence, trivially, $\Ker{(\Omega_n)}{\omega}=\{\omega\}$, or $\Ker{(\Omega_n)}{\omega}=G\neq\emptyset$. In the latter case, $\Ker{(\Omega_n)}{\omega}$ coincides with the connected component of
$
\bigcup_{n\in\N}\interior\big(\bigcap_{m\ge n}\Omega_{m}\big)
$
that contains~$\omega$. Here $\interior(\cdot)$ stands for the topological interior of a set.

As a matter of convenience, we also define the kernel of~$(\Omega_n)$ w.r.t. points ${\omega\in\C}$ that fail to belong to all but a finite number of $\Omega_n$'s. In such a case, we define ${\Ker{(\Omega_n)}{\omega}:=\{\omega\}}$ if there exists a sequence $(\omega_n)$ converging to~$\omega$ with $\omega_n\in\Omega_n$ for all ${n\in \N}$; otherwise, we put  ${\Ker{(\Omega_n)}{\omega}:=\emptyset}$.

The kernel of~$(\Omega_n)$ w.r.t.~$\omega$ is said to be \textsl{non-trivial} if it is different from~$\emptyset$ and~$\{\omega\}$. In such a case, $\Ker{(\Omega_n)}{\omega}$ is a domain in~$\C$ containing~$\omega$. Otherwise, i.e. if ${\Ker{(\Omega_n)}{\omega}\in\big\{\emptyset,\{\omega\}\big\}}$, we say that the kernel of~$(\Omega_n)$  w.r.t.~$\omega$ is \textsl{trivial}.

Note that for any subsequence $(\Omega_{n_k})$, $\Ker{(\Omega_{n_k})}{\omega}\supset \Ker{(\Omega_{n})}{\omega}$ and, in general, the inclusion can be strict.  A sequence $(\Omega_n)$ is said to \textsl{converge to its kernel $\Omega_*:=\Ker{(\Omega_n)}{\omega}$ w.r.t. a point~$\omega\in\C$,} if $\Omega_*\neq\emptyset$ and $\Ker{(\Omega_{n_k})}{\omega}=\Omega_*$ for every subsequence $(\Omega_{n_k})$.

\medskip

The above ``sequential'' concepts can be extended to continuous indexes in a natural way. Consider a family~$(\Omega_r)_{r>0}$ of domains in~$\C$ and let ${\omega\in\C}$. If for some $r_0>0$, a fixed neighbourhood of~$\omega$ is contained in~$\Omega_r$ whenever $r\geq r_0$, then $\Ker{(\Omega_r)}{\omega}$, the kernel of the family $(\Omega_r)$ w.r.t.~$\omega$, is defined as the connected component of
$$
\bigcup_{r>0}\,\interior\Big(\bigcap_{r'\ge r}\Omega_{r'}\Big)~=~\interior\Big(\bigcup_{r>0}\,\bigcap_{r'\ge r}\Omega_{r'}\Big)
$$
that contains~$\omega$. Otherwise, we put $\Ker{(\Omega_r)}{\omega}:=\{\omega\}$ or $\Ker{(\Omega_r)}{\omega}:=\emptyset$ depending on whether there exists a map $(0,+\infty)\ni r\mapsto \omega_r\in\C$ such that $\omega_r\in\Omega_r$ for all~$r>0$ and $\omega_r\to\omega$ as $r\to+\infty$.

The family $(\Omega_r)$ is said to \textsl{converge to its kernel $\Omega_*:=\Ker{(\Omega_r)}{\omega}$ w.r.t.~$\omega$} if $\emptyset\neq\Omega_*=\Ker{(\Omega_{r_n})}{\omega}$ for every sequence $(r_n)\subset(0,+\infty)$ tending to~$+\infty$.

\begin{remark}\label{RM_kernel-wellknown}
It follows easily from the definition, that if $K\subset\Ker{(\Omega_r)}{\omega}$ is a compact set, then $K\subset\Omega_r$ for all ${r>0}$  large enough. Conversely, if a domain~$U$ is contained in~$\Omega_r$ for all $r>0$ large enough, then $U\subset\Ker{(\Omega_r)}{\omega}$ for any $\omega\in U$. Analogous statements hold for kernels of sequences of domains.
\end{remark}

In the proof of our main result, Theorem~\ref{main}, we make use of the following statement, which is an easily corollary of Carath\'eodory's classical Kernel Convergence Theorem; see e.g. \cite[Theorem~1 in~\S{}II.5]{Goluzin}.

\begin{proposition}\label{PR_CarathConv}
Let $(g_n)$ be a sequence of conformal mappings of $\UD$ into~$\C$. If $(g_n)$ converges locally uniformly in~$\UD$ to some function~$g$, then $(g_n(\UD))$ converges to its kernel  w.r.t.~$\omega:=g(0)$. Moreover, $g(\UD)=\Ker{(g_n(\UD))}\omega$.

If the kernel $\Ker{(g_n(\UD))}{\omega}$ is non-trivial, then $g$ is conformal and on every compact set $K\subset g(\UD)$, the sequence $(g^{-1}_n)$ converges uniformly to~$g^{-1}$.
\end{proposition}

As a consequence of Remark~\ref{RM_Koenings-image}, in this paper, we will be especially interested in domains starlike at infinity. Simple ``model examples'' of such domains, relevant to the slope problem, are represented by angular sectors of the form
\[
S_p(\beta_1, \beta_2):=\big\{p+ite^{i\theta}\colon t>0,\; -\beta_2<\theta<\beta_1  \big\},
\]
where $p\in \C$ and $0\leq \beta_1,\beta_2\leq\pi$ with $\beta_1+\beta_2>0$.

\begin{remark}\label{Example} Clearly, when the above notions are applied to describing the limit behaviour of domains, much depends on the choice of the point~$\omega$ involved in the definition of the kernel.  Given a family $(\Omega_r)$ of domains and a sequence $(r_n)\subset(0,+\infty)$ tending to~$+\infty$, the limit behaviour of the sequence $(\Omega_{r_n})$ w.r.t. to some points $\omega\in\Complex$ can be similar to that of the whole family~$(\Omega_r)$, while for other choices of~$\omega$, $(\Omega_{r_n})$ and ~$(\Omega_r)$ can behave differently.  Consider the following example. Let $\beta\in(0,\pi]$,
 $$
 \Omega:=S_0(\pi/4,\beta)\,\Big\backslash\,\bigcup_{n=0}^{\infty}\,\Big\{u+iv\colon u=-2^n,~v\in[2^n,2^n(1+2^n)]\Big\},
 $$
 and define $\Omega_r:=\frac1r\Omega$ for all $r>0$.
 It can be checked that if $\omega\in S_0(0,\beta)$, then $\Ker{(\Omega_r)}\omega=S_0(0,\beta)$, and for all $\omega\in\Complex\setminus S_0(0,\beta)$ the kernel $\Ker{(\Omega_r)}\omega$ is trivial.  In particular, $S_0(0,\beta)$ is the unique non-trivial kernel of the family $(\Omega_r)$. Moreover, $(\Omega_r)$ converges to its kernel $S_0(0,\beta)$ w.r.t. any $\omega\in S_0(0,\beta)$.

It follows that  the sequence $(\Omega_{2^n})$ converges to its kernel $S_0(0,\beta)$ w.r.t. any ${\omega\in S_0(0,\beta)}$. However, $(\Omega_{2^n})$ has infinitely many other non-trivial kernels with respect to points in the left half-plane; namely,
$$
D_k:=\Ker{\{(\Omega_{2^n})}{\omega_k}=\{u+iv\colon {u\in(-2^k,-2^{k-1})}{,~} {v>-u}\},
$$
where $\omega_k:=(-\tfrac34+i)2^k,$ $k\in\mathbb{Z}$.
In fact, for any $k\in\mathbb{Z}$, the sequence $(\Omega_{2^n})$ converges to its kernel~$D_k$ w.r.t.~$\omega_k$.
Note also that in this example, $\Omega$ and hence all $\Omega_r$'s  are starlike at infinity.
 \end{remark}

For families $(\Omega_r)$ generated, as in the above remark, by scaling a given domain~$\Omega$,  the fact that the parameter $r$ takes all positive real values imposes strong restrictions on possible non-trivial kernels.

\begin{proposition}\label{LM_Omega-kernel}
	Let $\Omega\subset\C$ be a domain different from $\C$. Suppose that $(\Omega_r:=\frac1{r}\Omega)_{r>0}$ has a non-trivial kernel $\Omega_*:=\Ker{(\Omega_r)}{\omega}$ w.r.t. some point~$\omega\in\C$. Then the following assertions hold.
	\begin{itemize}
		\item[(A)] Either $\Omega_*$ coincides with $\C^*:=\C\setminus\{0\}$ or $\Omega_*$ is an angle with the vertex at the origin, i.e.
		$$
		\Omega_*=\lambda S_0(\beta_1,\beta_2),
		$$
		for some $\lambda\in\partial \D$ and some $0\leq \beta_1,\beta_2\leq\pi$ with $\beta_1+\beta_2>0$.\medskip
		
		\item[(B)] If, in addition, $\Omega$ is starlike at infinity, then
		$$\Omega_*= S_0(\beta_1,\beta_2),
		$$
		for some $0\leq \beta_1,\beta_2\leq\pi$ with $\beta_1+\beta_2>0$.
		Moreover, $\Omega_*$ is the only non-trivial kernel of $(\Omega_r)$, i.e. $\Ker{(\Omega_r)}{\omega'}\in\big\{\emptyset,\{\omega'\},\Omega_*\big\}$ for any $\omega'\in\C$.\medskip
	
		\item[(C)] Under hypothesis of~(B), suppose that  $\beta_1\beta_2\neq0$. Then for any $\omega'\in\C$ and any sequence $(r_n)\subset (0,+\infty)$ converging to $+\infty$ such that the kernel $\Ker{(\Omega_{r_n})}{\omega'}$ is non-trivial, we have
		$$
		\Ker{(\Omega_{r_n})}{\omega'}=\Ker{(\Omega_{r_n})}{\omega}\supset S_0(\beta_1,\beta_2).
		$$
\end{itemize}
\end{proposition}

\begin{proof}To prove~(A), we notice that $0\not\in \Omega_*$, because otherwise $\{z:|z|\leq\varepsilon\}\subset\Omega_r$ for some ${\varepsilon>0}$ and all~$r>0$ large enough and hence we would have $\Omega=\C$. Therefore, to prove~(A), it is enough to show that together with any $w\in\Omega_*$, the domain~$\Omega_*$ contains also the ray~$\{aw:a>0\}$. Suppose on the contrary that $w\in\Omega_*$ and that there exists $a>0$ such that $aw\in\partial\Omega_*$. Choose $\varepsilon>0$ so small that $\{\xi:|\xi-w|\leq\varepsilon\}\subset\Omega_*$. Then there exists $r>0$ such that for all $r'\ge r$, $\{\xi:|\xi-w|\leq\varepsilon\}\subset\Omega_{r'}=\frac1{r'}\Omega$ or equivalently, $\{a\xi:|a\xi-aw|\leq a\varepsilon\}\subset \frac{a}{r'}\Omega$. Hence  $U:=\{\xi:|\xi-aw|< a\varepsilon\}\subset\Omega_{t}$ for all~$t\ge r/a$. Since by construction the domain $U$ intersects $\Omega_*$, it follows that $U\subset\Omega_*$. We obtain a contradiction, which shows that the whole ray~$\{aw:a>0\}$ is contained in~$\Omega_*$.\medskip

\noindent
\textit{Proof of (B).} Let us now assume that $\Omega+it\subset\Omega$ for all ${t\ge0}$. Then the same property is possessed by $\Omega_r$ for any~$r>0$ and hence by $\Omega_*$. Since $0\not\in\Omega_*$, it follows that the ray $\{-it:t\ge0\}$ is contained in~$\C\setminus\Omega_*$. Hence $\Omega_*=S_0(\beta_1,\beta_2)$ for some $\beta_1,\beta_2\in[0,\pi]$. Clearly, $(\beta_1,\beta_2)\neq(0,0)$ because $\Omega_*\neq\emptyset$.

Take $\omega^\prime\in\C$ different from $\omega$ and assume  $\Ker{(\Omega_r)}{\omega'}$ is non-trivial. Note that $\Ker{(\Omega_r)}{\omega}$ and $\Ker{(\Omega_r)}{\omega'}$ either coincide or do not intersect.  Therefore, if both $\beta_1\neq0$ and $\beta_2\neq0$, then clearly,
$
\Ker{(\Omega_r)}{\omega^\prime}=\Ker{(\Omega_r)}{\omega}=S_0(\beta_1,\beta_2).
$

Let us show that the same conclusion holds also if $\beta_1\beta_2=0$. Otherwise, without loss of generality, we may assume that $\Ker{(\Omega_r)}{\omega}=S_0(\alpha_1,0)$ and $\Ker{(\Omega_r)}{\omega'}=S_0(0,\alpha_2)$ for some $\alpha_1,\alpha_2\in(0,\pi]$. Note these two sets are the only non-trivial kernels of $(\Omega_r)$. Hence
\begin{equation}\label{EQ_two-angles}
S_0(\alpha_1,0)\cup S_0(0,\alpha_2)~=~\interior\Big(\bigcup_{r>0}\,\bigcap_{r'\ge r}\Omega_{r'}\Big).
\end{equation}
In particular, $\Omega$ is neither contained in~$\{w:\Re w>0\}$ nor in~$\{w:\Re w<0\}$. Together with the property that $\Omega+it\subset\Omega$ for all $ {t\ge0}$, this means that  $\{i(t+a_0): t\ge0\}\subset\Omega$ for some $a_0\in\R$. If $a_0\le0$, then we immediately get a contradiction because in such a case $\{it:t\ge0\}\subset\Omega_r$ for any~$r>0$. If $a_0>0$, then for any $r>0$ we have $\{i(t+\frac{a_0}{r}):t\ge0\}\subset\Omega_r$. It follows that $\{it:t\ge0\}$ is contained in $\bigcup_{r>0}\bigcap_{r'\ge r}\Omega_{r'}$, which again contradicts~\eqref{EQ_two-angles}. This completes the proof of~(B).\medskip

\noindent \textit{Proof of~(C).} Fix now $\omega^\prime\in\C$ and a sequence $(r_n)\subset (0,+\infty)$ converging to $+\infty$ with the property that the kernel $\Ker{(\Omega_{r_n})}{\omega'}$ is non-trivial. Clearly $$
D:=\Ker{(\Omega_{r_n})}{\omega}\supset \Ker{(\Omega_{r})}{\omega}=S_0(\beta_1,\beta_2).
$$
As above, the domains $D$ and $D':=\Ker{(\Omega_{r_n})}{\omega'}$ either coincide or do not intersect; moreover, ${D'+it\subset D'}$ for all~$t\ge0$.
Since by hypothesis $\beta_1>0$ and $\beta_2>0$, it follows that  $D'=D$. The proof is now complete.\qed
\end{proof}
\begin{remark}
The condition   $\beta_1\beta_2\neq0$ in part~(C) of the above proposition is essential. Just consider the example in Remark~\ref{Example}.
\end{remark}

\begin{remark}\label{RM_shift-inv}
	Fix $w_0\in\C$ and let $\Omega\subset\C$ be a domain.  If  the family of domains $\widetilde\Omega_r:=\frac1r(\Omega-w_0)$ has a non-trivial kernel with respect to a point ${\omega\in\C}$, then also the family $\big(\Omega_r:=\frac1r\Omega\big)$ has the same kernel with respect to~$\omega$. A similar assertion holds for any sequence $(\Omega_{r_n})$ with ${r_n\to+\infty}$ as~${n\to+\infty}$. Indeed, every point  $$w\in\bigcup_{r>0}\interior\big(\bigcap_{r'\ge r}\widetilde\Omega_{r'}\big)$$ is contained in $\widetilde\Omega_{r'}$ along with some fixed neighbourhood for all~${r'>0}$ large enough. Hence $\{\xi:|\xi-w|<\varepsilon\}\subset\Omega_{r'}$ for some ${\varepsilon>0}$ and all ${r'>0}$ large enough. Taking into account that the same holds also with $\widetilde\Omega_{r'}$ and $\Omega_{r'}$ interchanged, we conclude that $$
	\bigcup_{r>0}\interior\big(\bigcap_{r'\ge r}\Omega_{r'}\big)=\bigcup_{r>0}\interior\big(\bigcap_{r'\ge r}\widetilde\Omega_{r'}\big).
	$$
\end{remark}

\section{Angular extent functions for domains starlike at infinity}
Recall that a domain  starlike at infinity  is a domain $\Omega$ of the complex plane such that $\Omega+it\subset \Omega$ for all $t\geq 0$.

\begin{definition} Let $\Omega$ be a domain  starlike at infinity.  A point $p\in\C$ is said to be a \textsl{natural point} associated with $\Omega$ if there exists $t_0\geq0$ such that $p+it\in \Omega$, for all $t>t_0$.
The set of all natural points associated with $\Omega$ will be denoted by $\textsf{NP}(\Omega)$.
\end{definition}

\begin{remark} The set $\textsf{NP}(\Omega)$ is always non-empty. Indeed, $\Omega\subset \textsf{NP}(\Omega)$. Moreover, $\textsf{NP}(\Omega)=\C$ if and only if $\{\Re z\colon z\in\Omega\}$ is unbounded both from above and from below, i.e. if and only if $\Omega$ is not contained in a half-plane bounded by a line parallel to the imaginary axis.
\end{remark}

\begin{remark}
Consider a non-elliptic one-parameter semigroup $(\varphi_t)$ in $\D$ with planar domain $\Omega$. It is known, see e.g. \cite[Theorem~9.3.5]{BCD-Book}, that the semigroup $(\varphi_t)$ is hyperbolic if and only if $\textsf{NP}(\Omega)$ is an open strip and that it is parabolic of positive hyperbolic step if and only if $\textsf{NP}(\Omega)$ is an open half-plane; finally, $(\varphi_t)$ is parabolic of zero hyperbolic step if and only if $\textsf{NP}(\Omega)=\C$.	
\end{remark}

\begin{remark}\label{RM_NP0} A useful subset of $\textsf{NP}(\Omega)$, which will be denoted by  $\textsf{NP}_0(\Omega)$, is formed by those points $p\in\partial\Omega$ such that $p+it\in \Omega$, for all $t>0$.
	
Note that $\textsf{NP}_0(\Omega)$ can be empty, but this happens only for a narrow class of domains of the form $\Omega=\vstrip(I):=\{x+iy\colon x\in I,\,y\in\R\}$,
where $I\subset\R$ is an open interval (bounded or unbounded or the whole~$\R$).

Moreover, $\textsf{NP}_0(\Omega)$ can be reduced to a unique point: consider, e.g. $\Omega:=\C\setminus\{-it:t\geq 0\}$.
\end{remark}

\begin{definition}\label{angular} Let $\Omega$ be a domain starlike at infinity. Fix $p\in \textsf{NP}(\Omega)$. For any $t>0$ such that $p+it\in\Omega$, we define the \textsl{(normalized) left angular extent} of $\Omega$ w.r.t.  $p$  by
\[
\alpha^-_{\Omega,p}(t):=\min\big\{\pi,\,\sup\{\alpha> 0:p+ite^{i\theta}\in \Omega \mathrm{\ for\  all\ } \theta\in[0,\alpha]\}\big\}\in(0,\pi].\hphantom{{}^-}
\]
Likewise, the \textsl{(normalized) right
angular extent} of $\Omega$ w.r.t. $p$ is defined as
\[
\alpha^+_{\Omega,p}(t):=\min\big\{\pi,\,\sup\{\alpha> 0:p+ite^{-i\theta}\in \Omega \mathrm{\ for\  all\ } \theta\in[0,\alpha]\}\big\}\in(0,\pi].
\]

The \textsl{natural domain of definition} for both functions $\alpha^-_{\Omega,p}$ and $\alpha^+_{\Omega,p}$ is the interval $\big(t_0(p),+\infty\big)$, where
\begin{equation}\label{EQ_naturalD}
t_0(p):=\inf\big\{t\ge0\colon p+it\in\Omega\big\}.
\end{equation}
\end{definition}

\begin{remark} Note that
$p\in \textsf{NP}_0\big(S_p(\beta_1,\beta_2)\big)$ and that for all $t>0$,
$$
\alpha^-_{S_p(\beta_1,\beta_2),p}(t)=\beta_1,\quad \alpha^+_{S_p(\beta_1,\beta_2),p}(t)=\beta_2.
$$
\end{remark}

\begin{remark} It is easy to see that for any domain $\Omega$ starlike at infinity and any $p\in \textsf{NP}(\Omega)$, both functions $\alpha^-_{\Omega,p}$ and $\alpha^+_{\Omega,p}$ are continuous from the right on their natural domain of definition. However, in general, neither $\alpha^-_{\Omega,p}$ nor $\alpha^+_{\Omega,p}$ are continuous from the left. For instance, for $\Omega:=\{z\colon \Im z>0\}\setminus\{-1+iy:y\in[0,1]\}$, we have
\[
\alpha_{\Omega,0}^-(t)=
\begin{cases}
\pi/2, & \text{if~}~t\in(0,1], \\
\arcsin(1/t), & \text{if~}~t\in(1,\sqrt{2}], \\
\pi/2, & \text{if~}~t\in(\sqrt{2},+\infty),
\end{cases}
\]
which is discontinuous from the left at $t=\sqrt{2}$.
\end{remark}

\begin{remark}
In general, for different points $p,q\in \textsf{NP}(\Omega)$, we have different functions $\alpha^+_{\Omega,p}$ and $\alpha^+_{\Omega,q}$. The same holds for the left angular extends. However, quite often the angular extends w.r.t. different points behave the same way in the limit as~${t\to+\infty}$.
\end{remark}

\begin{theorem}\label{kernel} Let $\Omega$ be a domain starlike at infinity and $p\in \textsf{NP}(\Omega)$. Let $t_0(p)$ be defined by~\eqref{EQ_naturalD}. Then the three conditions below are equivalent to each other.
\begin{itemize}
\item[(i)] The following limits exist,
$$
\alpha^-(p):=\lim_{t\to+\infty}\alpha^{-}_{\Omega,p}(t)\in[0,\pi], \qquad  \alpha^+(p):=\lim_{t\to+\infty}\alpha^{+}_{\Omega,p}(t)\in[0,\pi],
$$
and at least one of them is different from zero.\medskip

\item[(ii)] There are two increasing unbounded sequences $(t'_n),\, (t''_n)\subset\big(t_0(p),+\infty\big)$ with\\ ${t'_{n+1}/t'_n\to1}~$ and $~{t''_{n+1}/t''_n\to1}$ as ${n\to\infty}$ such that the following limits exist,
$$
\lim_{n\to\infty}\alpha^{-}_{\Omega,p}(t^\prime_n)\in[0,\pi], \qquad  \lim_{n\to\infty}\alpha^{+}_{\Omega,p}(t^{\prime\prime}_n)\in[0,\pi],
$$
and at least one of them is different from zero.\medskip

\item[(iii)] There exist $\omega\in\C$ with respect to which the family $\big(\Omega_r:=\frac1r\Omega\big)$ converges to a non-trivial kernel.
\end{itemize}
Moreover, if one and hence all of the above conditions are satisfied, then:
\begin{itemize}
\item[(a)] $\Ker{(\Omega_r)}{\omega'}=S_p\big(\alpha^-(p),\alpha^+(p)\big)$ for any $\omega'\in\C$ such that $\Ker{(\Omega_r)}{\omega'}$ is non-trivial, and
\item[(b)] for any $q\in\textsf{NP}(\Omega)$, $\alpha^{-}_{\Omega,q}(t)\to\alpha^{-}(p)$ and $\alpha^{+}_{\Omega,q}(t)\to\alpha^{+}(p)$ as $t\to+\infty$.
\end{itemize}
\end{theorem}

For the proof of this theorem we need one technical lemma.
\begin{lemma}\label{LM_close-radii}
	Under conditions of Theorem~\ref{kernel}, fix some $t>0$ and $\beta\in(0,\pi]$ such that  $L:=\{p+ite^{-i\theta}:0\le\theta<\beta\}\subset\Omega$. Then the following statements are true.
	\begin{itemize}
		\item[(A)] If $\,0<\beta\le \pi/2$, then for any $\varepsilon\in(0,\beta)$ there exists $\delta=\delta(\beta,\varepsilon)>0$ such that $\alpha_{\Omega,p}^+(x)\ge\beta-\varepsilon$ for all $x\in[t,t(1+\delta)]$.
		\medskip
		\item[(B)] If $\,\pi/2<\beta\le\pi$, then for any $\varepsilon\in(0,\beta)$ there exists $\delta=\delta(\beta,\varepsilon)\in(0,1)$ such that for each $x\in [t(1-\delta),t]$ at least one of the inequalities
		$$
		\alpha_{\Omega,p}^+(x)\le \pi-(\beta-\varepsilon)\quad \mathrm{\ or \ }\quad \alpha_{\Omega,p}^+(x)\ge \beta-\varepsilon
		$$
		holds.
	\end{itemize}
\end{lemma}
\begin{proof}
Note that the hypothesis implies that $t$ belongs to the natural domain of definition of~$\alpha_{\Omega,p}^-$ and~$\alpha_{\Omega,p}^+$.
For simplicity, we will assume, without loss of generality, that $p=0$.
		
Suppose first that $0<\beta\le \pi/2$.
Since $\Omega+it\subset\Omega$ for all ${t\ge0}$, together with the arc~$L$, the domain $\Omega$ contains the set
$$
     L_+:=\{z:\, 0\leq \Re z< t\sin (\beta), \ \Im z\geq \sqrt{t^{2}-(\Re z)^{2}}\}.
$$
Note that for any $x\ge t$, $$\{ixe^{-i\theta}\colon 0\le\theta<\theta(x)\}\subset L_+,\quad \theta(x):=\arcsin \frac{t\sin\beta}{x}. $$ The statement~(A) follows now easily.
	
	Similarly, if $\pi/2<\beta\le \pi$, then together with the arc $\{ite^{-i\theta}\colon\pi/2\le\theta<\beta\}$, the domain  $\Omega$ contains the set
	$$
	S:=\{z:|z|\leq t\}\,{\textstyle\bigcap}\,\{z:\Re z>t\sin\beta\}.
	$$
	In case $\beta=\pi$, it immediately follows that $\alpha_{p,\Omega}^+(x)=\pi$ for all $x\in(0,t]$. Suppose now that $\beta\in(\pi/2,\pi)$. Then for any $x\in(t\sin\beta,t]$, the arc $\{ixe^{-i\theta}:\theta(x)<\theta<\pi-\theta(x)\}$ is contained in~$S$, where $\theta(x)$ is defined as above. Therefore, for all such points $x$, either $\alpha_{\Omega,p}^+(x)\le\theta(x)$ or $\alpha_{\Omega,p}^+(x)\ge\pi-\theta(x)$. This implies statement~(B). \qed
\end{proof}

\begin{theopargself}
\begin{proof}[\!of Theorem~\ref{kernel}] To show the equivalence of the three conditions, it is sufficient to prove that (ii)~$\Longrightarrow$~(iii) and (iii)~$\Longrightarrow$~(i). Thanks to Remark~\ref{RM_shift-inv}, replacing the domain $\Omega$ with $\Omega-p$, we may assume that $p=0$.

\medskip\noindent\textit{Proof of (ii)~$\Longrightarrow$~(iii).}
Denote
$$
\alpha^-:=\lim_{n\to\infty}\alpha^{-}_{\Omega,0}(t^\prime_n), \qquad  \alpha^+:=\lim_{n\to\infty}\alpha^{+}_{\Omega,0}(t^{\prime\prime}_n).
$$

We are going to show that under condition~(ii) the following two claims hold.

\medskip\noindent\textit{Claim 1:} $\liminf_{t\to+\infty}\alpha_{\Omega,0}^-(t)\ge \alpha^-~$ and $~\,\liminf_{t\to+\infty}\alpha_{\Omega,0}^+(t)\ge \alpha^+$.

\medskip\noindent\textit{Claim 2:} for any sequence $\{r_n\}\subset(0,+\infty)$ tending to~$+\infty$ and any $w\in \partial S_0(\alpha^-,\alpha^+)$ there is a sequence $\{w_n\}$ converging to~$w$
and satisfying $w_n\in\C\setminus\Omega_{r_n}$ for all $n\in\N$.

\medskip
Claim~1 implies that for $\omega:=ie^{i(\alpha^--\alpha^+)/2}$ and for any sequence $(r_n)\subset (0,+\infty)$ with $r_n\to+\infty$ as $n\to+\infty$, the kernel $\Ker{(\Omega_{r_n}))}{\omega}$ contains $S_0(\alpha^-,\alpha^+)$. Likewise, from Claim~2, it follows that $\Ker{(\Omega_{r_n}))}{\omega}$ is contained in~$S_0(\alpha^-,\alpha^+)$.

\medskip\noindent\textit{Proof of Claim 1.} The proofs of the two inequalities in Claim~1 are similar. By this reason we prove only one of them, namely, the inequality for $\alpha_{\Omega,0}^+$. Suppose on the contrary that there exists $\alpha_*<\alpha^+$ and a sequence $x_n\to+\infty$ such that $\alpha_{\Omega,0}^+(x_n)\le \alpha_*$ for all~$n\in\N$. Consider two cases: $\alpha_*<\pi/2$ and $\alpha_*\in[\pi/2,\pi)$. Suppose first that ${\alpha_*<\pi/2}$. Then for all $n\in\N$ large enough, $\alpha_{\Omega,0}^+(t''_n)\ge \beta:=\min\{(\alpha_*+\alpha^+)/2,\pi/2\}$. Note that $\beta>\alpha_*$. Hence applying Lemma~\ref{LM_close-radii}\,(A) for $t:=t''_n$, we conclude that there exist $n_0\in\N$ and $\delta>0$ such that for all $x\in J:=\bigcup_{n\ge n_0}[t''_n,t''_n(1+\delta)]$ we have $\alpha_{\Omega,0}^+(x)>\alpha_*$. Since $t''_{n+1}/t''_n\to1$ as $n\to+\infty$,  $[a_0,+\infty)\subset J$ for some $a_0>0$, and in particular, $x_n\in J$ for all sufficiently large ${n\in\N}$. This contradicts the assumption that $\alpha_{\Omega,0}^+(x_n)\le \alpha_*$ for all~$n\in\N$.

Let us now obtain a contradiction in the case $\alpha_*\in[\pi/2,\pi)$. Since  $\alpha_*<\alpha^+$, in this case $\alpha^+>\pi/2$. It follows that for all $n\in\N$ large enough, $\alpha_{\Omega,0}^+(t''_n)\ge\pi/2$. Therefore, similarly to the first case, applying Lemma~\ref{LM_close-radii}\,(A) for $t:=t''_n$ and $\beta:=\pi/2$, one can show that there exists $a_1>0$ such that for all  $x\ge a_1$, we have $\alpha_{\Omega,0}^+(x)>\pi-\beta'$ with $\beta':=(\alpha_*+\alpha^+)/2$.  Moreover, taking into account that $\alpha_{\Omega,0}^+(t''_n)\ge\beta'$ for all $n\in\N$ large enough and  applying Lemma~\ref{LM_close-radii}\,(B) for $t:=t''_n$ and $\beta:=\beta'$, we see that there exists $n_0\in\N$  and $\delta\in(0,1)$ such that for all $x\in J':=\bigcup_{n\ge n_0}[t''_n(1-\delta),t''_n]$ we have $\alpha_{\Omega,0}^+(x)>\alpha_*$. This is again in contradiction with the assumption that $\alpha_{\Omega,0}^+(x_n)\le \alpha_*$ for all~$n\in\N$. Claim~1 is now proved.

\medskip\noindent\textit{Proof of Claim 2.} Clearly, it is sufficient to prove the claim for $w\neq0$. Then $w=i\rho e^{-i\theta}$ for some ${\rho>0}$ and $\theta\in\{-\alpha^-,\,\alpha^+\}$. We provide the proof only for the case ${\theta=\alpha^+}$, because the case ${\theta=-\alpha^-}$ is similar. Suppose first that $\alpha^+<\pi$. Then omitting if necessary a finite number of terms in~$(t''_n)$, we may suppose that $it_n''e^{-i\alpha_{\Omega,0}^+(t''_n)}\in\C\setminus\Omega$ for all~$n\in\N$.  For $x>0$, denote by $t(x)$ the element of the sequence~$(t''_n)$ for which $|\log(t''_n/x)|$ attains its minimal value $q(x)$. If the minimal value is attained for two different elements, we choose one of them, for example the smaller one. Clearly $\zeta_n:=it(\rho r_n)e^{-i\alpha_{\Omega,0}^+(t(\rho r_n))}\in\C\setminus\Omega$ for each~$n\in\N$. It follows that $w_n:=\zeta_n/r_n\in\C\setminus\Omega_{r_n}$. Moreover, $|\Re\,\log(w_n/w)|=q(\rho r_n)\to0$ because $t''_{n+1}/t''_n\to1$ as $n\to+\infty$, and $|\Im\,\log(w_n/w)|=|\alpha^+_{\Omega,0}(t(\rho r_n))-\alpha^+|\to0$ because $\alpha^+_{\Omega,0}(t''_n)\to\alpha^+$ as $n\to+\infty$. Thus $w_n\to w$ as $n\to+\infty$.

Suppose now that $\alpha^+=\pi$. The above argument works also in this case if ${0\not\in\Omega}$, but if ${0\in\Omega}$, then it might happen that $\alpha_{\Omega,0}^+(t)=\pi$ and $it e^{-i\alpha_{\Omega,0}^+(t)}=-it\in\Omega$ for all $t>0$. Hence we have to modify the above argument. To this end, fix some $w_*\in\C\setminus\Omega$. Then $\Gamma:=\{w_*-it:t\ge0\}\subset \C\setminus\Omega$. Removing a finite number of terms in~$(r_n)$,  we may suppose that $\rho r_n>|w_*|$ for all~$n\in\N$. Define $w_n$ to be the unique point of intersection $\Gamma\bigcap\,\{z:|z|=\rho r_n\}$ lying in the lower half-plane. Then clearly $\zeta_n:=w_n/r_n\to-i\rho=w$ as $n\to+\infty$. By construction, $\zeta_n\in\C\setminus\Omega_{r_n}$ for all ${n\in\N}$. Now Claim~2 and hence the implication (ii)~$\Longrightarrow$~(iii) is proved.

\medskip\noindent\textit{Proof of (iii)~$\Longrightarrow$~(i).}  By Proposition~\ref{LM_Omega-kernel}, $\Omega_*:=\Ker{(\Omega_r)}{0}=S_0(\beta_1,\beta_2)$ for some $\beta_1,\beta_2\in[0,\pi]$ with $\beta_1+\beta_2>0$. Since any compact subset of~$\Omega_*$ is contained in~$\Omega_r$ for all ${r>0}$ large enough, it follows that
$$
\liminf_{t\to+\infty}\alpha^{-}_{\Omega,0}(t)\ge\beta_1 \quad \text{\ and \ }\quad  \liminf_{t\to+\infty}\alpha^{+}_{\Omega,0}(t)\ge\beta_2.
$$
It remains to show that
$$
 \limsup_{t\to+\infty}\alpha^{-}_{\Omega,0}(t)\le\beta_1\quad \text{\ and \ }\quad \limsup_{t\to+\infty}\alpha^{+}_{\Omega,0}(t)\le\beta_2.
$$
We are going to prove only the latter claim, since the proof of the other one is very similar. Thus suppose that $\limsup_{t\to+\infty}\alpha^{+}_{\Omega,0}(t_n)>\beta_2$. Then there exists $\beta_*>\beta_2$ and an unbounded strictly increasing sequence $(t_n)\subset(0,+\infty)$ such that ${\alpha^{+}_{\Omega,0}(t_n)\ge\beta_*}$ for all ${n\in\Natural}$. If $\beta_2<\pi/2$, then using assertion~(A) of Lemma~\ref{LM_close-radii} with $\beta:=\min\{\beta_*,\pi/2\}$, we see that there is $\mu\in(0,1)$ and a sequence of intervals $I_n:={[r_n(1-\mu),r_n(1+\mu)]}$ with ${r_n\to+\infty}$ as ${n\to+\infty}$ such that $\alpha^{+}_{\Omega,0}(x)\ge \beta':=(\beta+\beta_2)/2$ for any ${x\in\bigcup_{n\in\N}I_n}$. The same conclusion, but with $\beta':=(\beta_*+\beta_2)/2$, can be obtained in the case ${\beta_2\ge \pi/2}$ by applying assertion (B) of Lemma~\ref{LM_close-radii} with $\beta$ replaced by~$\beta_*$, if we recall that $\liminf_{t\to+\infty}\alpha^{+}_{\Omega,0}(t)\ge\beta_2\ge\pi/2$ and hence $\alpha^{+}_{\Omega,0}(t)>\pi-\beta_*$ for all $t>0$ large enough.

It follows that in both cases the set $U:=\{i\rho e^{-i\theta}:1-\mu\le\rho\le1+\mu,~ 0\le\theta<\beta'\}$ is contained in $\Omega_{r_n}$  for any $n\in\N$. Hence $\Ker{(\Omega_{r_n})}{0}\supset \Omega_*\bigcup\, U$. However, by~(iii) and the definition of convergence to the kernel, $\Ker{(\Omega_{r_n})}{0}=\Omega_*$. Since $U\not\subset\Omega_*$, we have obtained a contradiction, which means that indeed $\limsup_{t\to+\infty}\alpha^{+}_{\Omega,0}(t)\le\beta_2$.\medskip

\medskip\noindent\textit{Proof of (a) and (b).}
We have already seen that if conditions (i)\,--\,(iii) hold, then the equality
$
\Ker{(\Omega_r)}{\omega}=S_0\big(\alpha^-(p),\alpha^+(p)\big)
$
takes place for at least one point~$\omega\in\C$.
By Proposition~\ref{LM_Omega-kernel}, this implies assertion~(a).

Now (b) follows from (a) and the fact that condition~(iii) does not depend on the choice of the point~$p$ and take into account assertion~(a).\qed	
\end{proof}
\end{theopargself}
\begin{remark}\label{RM_NP-not-C} If $\Omega$ is a domain starlike at infinity with $\textsf{NP}(\Omega)\neq \C$, then one  of the following three mutually exclusive possibilities holds.
\begin{enumerate}
\item $\Omega\subset\vstrip(I):=\{x+iy\colon x\in I,~y\in\Real\}$ for a suitable \textit{bounded} interval ${I\subset\Real}$. In this case, for any $p\in \textsf{NP}(\Omega)$, we have $\lim_{t\to +\infty}\alpha^-_{\Omega,p}(t)=\lim_{t\to +\infty}\alpha^+_{\alpha,p}(t)=0.$
\item $\Omega\subset\vstrip(I)$ for a suitable interval $I$ of the form~$(a,+\infty)$, $a\in\Real$, but not for any bounded interval~$I$.  In this case, clearly,
$\lim_{t\to +\infty}\alpha^-_{\Omega,p}(t)=0$ for any $p\in \textsf{NP}(\Omega)$.
\item $\Omega\subset\vstrip(I)$ for a suitable interval $I$ of the form~$(-\infty,b)$, $b\in\Real$, but not for any bounded interval~$I$.  In this case,
$\lim_{t\to +\infty}\alpha^+_{\Omega,p}(t)=0$ for any $p\in \textsf{NP}(\Omega)$.
\end{enumerate}	
\end{remark}

\section{Boundary distance functions for domains starlike at infinity}
\begin{theopargself}
\begin{definition}[\protect{\cite{BCDGZ}}]~\,
Let $\Omega$ be a domain starlike at infinity, $p\in\C$ and $t>0$. The \textsl{(normalized) left distance} of $\Omega$ w.r.t. $p$ is defined by	
\[
\delta_{\Omega,p}^-(t):=\min\big\{t,~\inf\{|z-(p+it)|\colon \Re z\leq \Re p,~ z\in \C\setminus \Omega\}\big\}\in[0,t].
\]
Likewise, the \textsl{(normalized) right distance} of $\Omega$ w.r.t. $p$ is defined by
\[
\delta_{\Omega,p}^+(t):=\min\big\{t,~\inf\{|z-(p+it)|\colon \Re z\geq \Re p,~ z\in \C\setminus \Omega\}\}\in[0,t].
\]
\end{definition}
\end{theopargself}

\begin{remark} Note that if $p\in \textsf{NP}_0(\Omega)$, then, for all $t>0$,
$$
\delta_{\Omega,p}^-(t)=\inf\{|z-(p+it)|\colon \Re z\leq \Re p,~ z\in \partial \Omega\}\mathrlap{\quad\text{and}}
$$
$$
\delta_{\Omega,p}^+(t)=\inf\{|z-(p+it)|\colon \Re z\geq \Re p,~ z\in \partial \Omega\}.
$$	
\end{remark}

\begin{remark}
In contrast to the angular extents $\alpha_{\Omega,p}^\pm$, for any domain~$\Omega$ starlike at infinity and any ${p\in\C}$, the functions $\delta^\pm_{\Omega,p}$ are continuous and non-decreasing on the whole interval $(0,+\infty)$.
\end{remark}

\begin{theopargself}
\begin{theorem}[\protect{\cite[Lemma 3.6]{BCDGZ}}] \label{Lem:same-delta}
Let $\Omega$ be a domain starlike at infinity. Then for any ${p, q\in \C}$, there exist constants $c_2>c_1>0$  such that for all $t>0$,
\begin{equation*}
c_1 \delta^{-}_{\Omega,p}(t)\leq \delta^{-}_{\Omega,q}(t)\leq c_2\delta^{-}_{\Omega,p}(t)\mathrlap{\quad\text{and}}
\end{equation*}
\begin{equation*}
c_1 \delta^{+}_{\Omega,p}(t)\leq \delta^{+}_{\Omega,q}(t)\leq c_2\delta^{+}_{\Omega,p}(t).
\end{equation*}
\end{theorem}
\end{theopargself}

The following result obtained in~\cite{BCDGZ} establishes a strong relationship between the slopes of the trajectories of a one-parameter semigroup at its DW-point and the limit behaviour of the boundary distance functions of the corresponding planar domain.
\begin{theopargself}
\begin{theorem}[\protect{\cite[Theorem 1.1]{BCDGZ}}] \label{maintheoremBCDGZ}
	Let $(\varphi_t)$ be a non-elliptic semigroup in $\D$ with the DW-point $\tau\in\partial \D$ and Koenigs function $h$ and let $\Omega:=h(\D)$. Fix any sequence $(t_n)\subset(0,+\infty)$ tending to $+\infty$. Then:
	\begin{enumerate}
		\item[(A)] The sequence $(\varphi_{t_n}(z))$ converges non-tangentially to $\tau$ for some (and hence all)  ${z\in \D}$ if and only if for some (and hence all) ${p\in \Omega}$  there exist constants ${c_2>c_1>0}$ such that for all ${n\in \N}$,
		\[
		c_1\delta^+_{\Omega,p}(t_n)\leq \delta^-_{\Omega,p}(t_n)\leq c_2\delta^+_{\Omega,p}(t_n).
		\]
		\item[(B)]  $\lim_{n\to +\infty}\Arg(1-\overline{\tau}\varphi_{t_n}(z))={\pi}/{2}$ (in particular,  $(\varphi_{t_n}(z))$ converges tangentially to $\tau$ as ${n\to +\infty}$) for some (and hence all) $z\in \D$ if and only if for some (and hence all) $p\in \Omega$,
		\[
		\lim_{n\to +\infty}\frac{\delta^+_{\Omega,p}(t_n)}{\delta^-_{\Omega,p}(t_n)}=0,
		\]
        \item[(C)] $\lim_{n\to +\infty}\Arg(1-\overline{\tau}\varphi_{t_n}(z))=-{\pi}/{2}$ (in particular, $(\varphi_{t_n}(z))$ converges tangentially to $\tau$ as ${n\to +\infty}$) for some (and hence all) ${z\in \D}$ if and only if for some (and hence all) $p\in \Omega$,
		\[
		\lim_{n\to +\infty}\frac{\delta^+_{\Omega,p}(t_n)}{\delta^-_{\Omega,p}(t_n)}=+\infty.
		\]
	\end{enumerate}
\end{theorem}
\end{theopargself}

As one might expect, for domains starlike at infinity, angular extent functions and boundary distance functions are closely related, see Proposition~\ref{angular-distance} below. At the same time, it is worth to mention that these two characteristics are not asymptotically equivalent, as demonstrated by Example~\ref{EX_Santiago} in the last section. Therefore, the information on the geometry of the planar domain near~$\infty$ provided by the angular extents is not identical to that contained in the boundary distance functions.

\begin{proposition}\label{angular-distance}
Let $\Omega$ be a domain starlike at infinity. Fix some $p\in \textsf{NP}(\Omega)$ and let $(t_0,+\infty)$ be the natural domain of definition of $\alpha_{p,\Omega}^-$ and $\alpha_{p,\Omega}^+$. Then for all $t>t_0$,
\begin{enumerate}
\item[(a)] $\delta^+_{\Omega,p}(t) \leq2t\sin\big(\tfrac12\alpha^+_{\Omega,p}(t)\big)<t\alpha^+_{\Omega,p}(t)$;\medskip
\item[(b)] $\tfrac1\pi t\alpha^+_{\Omega,p}(t) \le t\sin\,\min\{\tfrac\pi2,\,\alpha^+_{\Omega,p}(t)\}\leq \delta^+_{\Omega,p}(2t)$;\medskip
\item[(c)] $\delta^-_{\Omega,p}(t) \leq2t\sin\big(\tfrac12\alpha^-_{\Omega,p}(t)\big)<t\alpha^-_{\Omega,p}(t)$;\medskip
\item[(d)] $\tfrac1\pi t\alpha^-_{\Omega,p}(t) \le t\sin\,\min\{\tfrac\pi2,\,\alpha^-_{\Omega,p}(t)\}\leq \delta^-_{\Omega,p}(2t)$.
\end{enumerate}
\end{proposition}
\begin{proof} Clearly, it is sufficient to prove (a) and~(b). The proof of (c) and~(d) is similar.
Without lose of generality, we assume that $p=0$. Moreover, to simplify the notation, for $t>t_0$, we will write $ \alpha^+(t):= \alpha^+_{\Omega,0}(t)$ and $\delta^+(t) :=\delta^+_{\Omega,0}(t)$.
	
\noindent\textit{Proof of (a)}. By the very definition, $\delta^+(t)\le t$. Hence (a) holds trivially if ${\alpha^+(t)=\pi}$. Therefore, we may suppose that ${\alpha^+(t)<\pi}$. In such a case, $w_0:=ite^{-\alpha^+(t)}\in\partial\Omega$ and we immediately get
$$
\delta^+(t)\le\big|it-w_0\big|=2t\sin\big(\tfrac12\alpha^+(t)\big).
$$

\noindent\textit{Proof of (b)}. Let $\beta:=\min\{\tfrac\pi2,\,\alpha^+(t)\}$. Note that the arc $\{ite^{-i\theta}: 0\leq \theta<\beta\}$ is a subset of $\Omega$. Since $\Omega$ is starlike at infinity, it follows that
$$
\left\{z\colon 0\leq \Re z< t\sin\beta,~ \Im z\ge t\right\}\subset \Omega.
$$
Therefore, $\delta^+(2t)=\min\big\{2t,\ \inf\{|z-2it|: z\in \C\setminus\Omega,\, \Re z\geq 0\}\big\}\geq   t\sin\beta.$\qed
\end{proof}

\section{Main results}
In this section we prove our main results, which establish relationships between the trajectory slopes at the DW-point and the asymptotic behaviour of the angular extents in the planar domain of the semigroup for ${t\to+\infty}$.

As we mentioned in the introduction, essentially the slope problem has been solved for hyperbolic semigroups and for parabolic semigroups of positive hyperbolic step. Therefore, we might strict our attention to parabolic semigroups of \textit{zero hyperbolic step}. At the same time, our methods do not require this assumption. That is why we will keep supposing only that the semigroup is non-elliptic.

We start with two corollaries of the main theorem from~\cite{BCDGZ}.
\begin{proposition}\label{prop-pendientes} Let $(\varphi_t)$ be a non-elliptic semigroup in $\D$  with DW-point $\tau\in\partial\D$, Koenigs function $h$, and planar domain $\Omega:=h(\D)$. Fix some $p\in \textsf{NP}(\Omega)$. Suppose that
\begin{equation}\label{EQ_tt2}
\liminf_{t\to+\infty}\frac{\alpha^-_{\Omega,p}(t/2)}{\alpha^-_{\Omega,p}(t)}>0\quad \mathrm{and}\quad \liminf_{t\to+\infty}\frac{\alpha^+_{\Omega,p}(t/2)}{\alpha^+_{\Omega,p}(t)}>0.
\end{equation}
Then the following are equivalent:
\begin{enumerate}
\item[(i)] The trajectory $t\mapsto\varphi_t(z)$ converges non-tangentially to $\tau$ as $t\to+\infty$ for some (and hence all) $z\in \D$.
\item[(ii)]  There exist $T>0$ and $0<C_1<C_2$ such that for all $t>T$,
\begin{equation}\label{EQ_C1C2-alpha}
C_1 \alpha^+_{\Omega,p}(t)\leq \alpha^-_{\Omega,p}(t)\leq C_2 \alpha^+_{\Omega,p}(t).
\end{equation}
\end{enumerate}	
\end{proposition}

\begin{proof}
Let $(t_0,+\infty)$ be the natural domain of definition of $\alpha_{\Omega,p}^-$ and $\alpha_{\Omega,p}^+$. To simplify the notation, for $t>t_0$, we will write $ \alpha^+(t):= \alpha^+_{\Omega,p}(t)$ and $\delta^+(t) :=\delta^+_{\Omega,p}(t)$. According to~\eqref{EQ_tt2}, there exists a constant $\varepsilon>0$ such that
\begin{equation}\label{EQ_tt2bis}
\alpha^-(t/2)\ge \varepsilon \alpha^-(t)\quad\text{and}\quad \alpha^+(t/2)\ge \varepsilon \alpha^+(t)
\end{equation}
for all $t>t_0$ large enough.

Suppose that (ii) holds. Then combining \eqref{EQ_C1C2-alpha}, \eqref{EQ_tt2bis}, and Proposition~\ref{angular-distance}, for all~$t>t_0$ large enough we obtain:
\begin{align*}
  &\delta^-(t)\ge t\pi^{-1}\alpha^-(t/2)\ge t\pi^{-1}C_1\alpha^+(t/2)\ge t\pi^{-1}C_1\varepsilon\alpha^+(t)\ge \pi^{-1} C_1\varepsilon\delta^+(t)\\
 \text{and}\quad &\delta^-(t)<t\alpha^-(t)\le t C_2\alpha^+(t)\le t\varepsilon^{-1} C_2\alpha^+(t/2)\le \varepsilon^{-1} C_2\pi\delta^+(t).
\end{align*}
Therefore, (i) holds by Theorem~\ref{maintheoremBCDGZ}\,(A).

Suppose now that condition~(i) is satisfied. Then applying again Theorem~\ref{maintheoremBCDGZ}\,(A), we see that for there exists $T\ge t_0$ and constants $c_2>c_1>0$ such that
$$
c_1\delta^+(t)\leq \delta^-(t)\leq c_2\delta^+(t)\quad \text{for all~$t>T$}.
$$
Combining these inequalites with \eqref{EQ_tt2} and Proposition~\ref{angular-distance}, we find that
\begin{align*}
& t\alpha^-(t)>\delta^-(t)\ge c_1\delta^+(t)\ge c_1\pi^{-1}t\alpha^+(t/2)\ge c_1\pi^{-1}t\varepsilon\alpha^+(t/2)\\
\text{and}\quad & t\alpha^-(t)\le t\varepsilon^{-1}\alpha^-(t/2)\le \varepsilon^{-1}\pi\delta^-(t)\le \varepsilon^{-1}\pi c_2\delta^+(t)<\varepsilon^{-1}\pi c_2 t\alpha^+(t)
\end{align*}
for all~$t>T$. It follows that~\eqref{EQ_C1C2-alpha} holds with $C_1:=c_1\pi^{-1}\varepsilon$ and $C_2:=\varepsilon^{-1}\pi c_2$.\qed
\end{proof}

\begin{remark}
Example~\ref{EX_Manolo} in the next section shows that condition~\eqref{EQ_tt2} in Proposition~\ref{prop-pendientes} is essential.
\end{remark}

\begin{proposition}\label{PR_second} Let $(\varphi_t)$ be a non-elliptic semigroup in $\D$  with DW-point $\tau\in\partial\D$, Koenigs function $h$, and planar domain $\Omega:=h(\D)$. Fix some $p\in \textsf{NP}(\Omega)$ and denote
$$
\alpha_p^-:=\liminf_{t\to+\infty}\alpha^-_{\Omega,p}(t)\in[0,\pi]\quad \mathrm{and}\quad \alpha_p^+:=\liminf_{t\to+\infty}\alpha^+_{\Omega,p}(t)\in[0,\pi].
$$
Then the following assertions hold:
\begin{enumerate}
\item[(A)] If $\alpha_p^->0$ and $\alpha_p^+>0$, then  the trajectories $t\mapsto\varphi_t(z)$ converge to~$\tau$ non-tangentially for all $z\in \D$.\smallskip
\item[(B)] If $\alpha_p^->0$ but $\alpha_{\Omega,p}^+(t)\to0$ as $t\to+\infty$, then $\mathrm{Slope}[t\mapsto \varphi_t(z),\tau]=\{\pi/2\}$ for every~${z\in\D}$. In particular, the trajectories $t\mapsto\varphi_t(z)$ converge to~$\tau$ tangentially  for all $z\in\D$.\smallskip
\item[(C)] If $\alpha_p^+>0$ but $\alpha_{\Omega,p}^-(t)\to0$ as $t\to+\infty$, then $\mathrm{Slope}[t\mapsto \varphi_t(z),\tau]=\{-\pi/2\}$ for every~${z\in\D}$. In particular, the trajectories $t\mapsto\varphi_t(z)$ converge to~$\tau$ tangentially  for all $z\in\D$.
\end{enumerate}
\end{proposition}
\begin{proof} Assertion~(A) is a corollary of Proposition~\ref{prop-pendientes}. Indeed, using a simple observation that
$$
 0< \frac{\alpha^-_{\Omega,p}(t)}{2\pi},\,\frac{\alpha^+_{\Omega,p}(t)}{2\pi}\le\frac12,
$$
we see that under the hypothesis of~(A), for all $t>0$ large enough we have
$$
\alpha^-_{\Omega,p}(t/2)\geq \frac{\alpha^-_{\Omega,p}(t)}{2\pi} \alpha^-_p,\quad \alpha^+_{\Omega,p}(t/2)\geq \frac{\alpha^+_{\Omega,p}(t)}{2\pi} \alpha^+_p,
$$
which implies~\eqref{EQ_tt2}, and
$$
\alpha^-_{\Omega,p}(t)\geq \frac{\alpha^+_{\Omega,p}(t)}{2\pi}{\alpha^-_p},\quad
\alpha^+_{\Omega,p}(t)\geq \frac{\alpha^-_{\Omega,p}(t)}{2\pi}{\alpha^+_p},
$$
which implies~\eqref{EQ_C1C2-alpha} for suitable $C_2>C_1>0$.

\smallskip\noindent
\textit{Proof of~(B)}. Since $\alpha^-(p)>0$, by Proposition~\ref{angular-distance}\,(d), $\delta^-_{\Omega,p}(t)>\varepsilon t$ for some ${\varepsilon>0}$ and all ${t>0}$ large enough. On the other hand, $\alpha^+_{\Omega,p}(t)\to0$ as ${t\to+\infty}$ and hence by Proposition~\ref{angular-distance}\,(a), $\delta^+_{\Omega,p}(t)/t\to0$ as ${t\to+\infty}$. Therefore, by Theorem~\ref{maintheoremBCDGZ}\,(B),\linebreak $\mathrm{Slope}[t\mapsto \varphi_t(z),\tau]=\{\pi/2\}$ and we are done.

\medskip
\noindent
The \textit{proof of~(C)} is analogous to that of~(B). Therefore, we may omit it.\qed
\end{proof}

As we have already mentioned, even asymptotically, the angular extends $\alpha^\pm_{\Omega,p}$ are not equivalent to the distance functions $\delta^\pm_{\Omega,p}$ used in~\cite{BCDGZ} (see Example~\ref{EX_Santiago}). In fact, we are able to establish the following result, which seems to have no analogues in terms of the distance functions.

\begin{theorem}\label{main} Let $(\varphi_t)$ be a non-elliptic semigroup in $\D$ with DW-point $\tau\in\partial\D$ and Koenigs function $h$. Let $\Omega:=h(\D)$. If for some $p\in \textsf{NP}(\Omega)$,
\[
\alpha^-(p):=\lim_{t\to +\infty}\alpha^-_{\Omega,p}(t)>0\quad \mathrm{and}\quad
\alpha^+(p):=\lim_{t\to +\infty}\alpha^+_{\Omega,p}(t)>0,
\]
then  $\mathrm{Slope}[t\mapsto \varphi_t(z),\tau]=\big\{\eta\tfrac\pi2\big\}$ for any $z\in\D$, where
$$
\eta:=\frac{\alpha^-(p)-\alpha^+(p)}{\alpha^-(p)+\alpha^+(p)}.
$$
In particular, every trajectory $t\mapsto\varphi_t(z)$ converges to $\tau$ non-tangentially and with a definite slope. 	
\end{theorem}
\begin{remark}
The conclusion in the above theorem, except for the non-tangential character of the trajectory convergence, remains valid when one of the limits $\alpha^\pm(p)$ is positive and the other is zero. This fact is a direct consequence of Proposition~\ref{PR_second}, but it can also be established independently using a technique similar to that we employ in the proof of Theorem~\ref{main}.
\end{remark}

\begin{remark}\label{RM_p-q}  It is worth pointing out that according to Theorem~\ref{kernel}, under the hypothesis of Theorem~\ref{main} for any other point $q\in \textsf{NP}(\Omega)$, the following limits exist
\[
\lim_{t\to +\infty}\alpha^-_{\Omega,q}(t),\quad
\lim_{t\to +\infty}\alpha^+_{\Omega,q}(t)
\]
and they coincide with $\alpha^-(p)$ and $\alpha^+(p)$, respectively.
\end{remark}

For the proof of  Theorem~\ref{main}, the following easy fact will be used:

\begin{remark}\label{carac_slope} Let $(\varphi_t)$ be a semigroup in $\D$ and $\tau\in\partial\D$. For a sequence $(t_n)\subset[0,+\infty)$ converging to $+\infty$, the following assertions are equivalent:
	\begin{enumerate}
		\item[(i)] There exists the limit $\theta :=\lim\limits_{n\to +\infty}\mathrm{Arg}\left( 1-\overline{\tau}\varphi_{t_n}(z)
		\right)\in \left[ -\frac{\pi }{2},\frac{\pi }{2}\right]$.
		\smallskip
		\item[(ii)]  There exists the limit $m:=\lim\limits_{n\to +\infty}\dfrac{1-\overline{\tau}\varphi_{t_n}(z)}{\left| 1-	\overline{\tau}\varphi_{t_n}(z)\right|}\in\partial\D$.
		\smallskip
		\item[(iii)]  There exists the limit $\mu :=\lim\limits_{n\to +\infty}\dfrac{\Im (\overline{\tau}\varphi_{t_n}(z))}{1-\Re (\overline{\tau}\varphi_{t_n}(z))}\in\left[-\infty,+\infty\right].$
	\end{enumerate}
	Moreover, if one and hence all of the above hold, then $e^{i\theta }=m$ and $
	\mu =-\tan \theta.$
\end{remark}

\begin{theopargself}
\begin{proof}[\textit{of Theorem~\ref{main}}] Without loss of generality we will assume that ${\tau=1}$.
Since $\alpha^+(p)>0$ and $\alpha^-(p)>0$, we have $\mathsf{NP}(\Omega)=\C$, see Remark~\ref{RM_NP-not-C}. In particular, it follows that $(\varphi_t)$ is a parabolic semigroup of zero hyperbolic step and that
for any constant ${c\in\Complex}$, $h+c$ is also a Koenigs function for $(\varphi_t)$. Therefore, bearing in mind Remarks~\ref{RM_NP0} and~\ref{RM_p-q}, we may also assume that $p=0\in\mathsf{NP}_0(\Omega)$. Then the natural domain of definition of $\alpha^-_{\Omega,0}$ and $\alpha^+_{\Omega,0}$ is $(0,+\infty)$.
Denote $$\alpha^-:=\lim_{t\to +\infty}\alpha^-_{\Omega,0}(t),\quad
\alpha^+:=\lim_{t\to +\infty}\alpha^+_{\Omega,0}(t),\quad \text{and~}~ \eta:=\dfrac{\alpha^- -\alpha^+}{\alpha^- +\alpha^+}.$$

Take $z_0\in\D$ such that $h(z_0)=i$. Thanks to Abel's equation~\eqref{EQ_Abel}, Theorem~\ref{indepen}, and Remark~\ref{carac_slope}, it is sufficient to show that
\begin{equation}\label{crucial}
\lim_{t\to+\infty}\frac{1-h^{-1}(it)}{|1-h^{-1}(it)|}=\exp\big(i\eta\tfrac\pi2\big).
\end{equation}

Since $0\in\partial \Omega$ and $h$ is univalent, the function $H(z):=-{1}/{h(z)}$ is holomorphic and univalent in~$\D$. Note that $\{it:t>0\}\subset H(\D)$, $\{-it:t\geq0\}\subset \C\setminus H(\D)$
 and $H^{-1}(w)=h^{-1}(-1/w)$ for all ${w\in H(\D)}$.

Consider the Jordan arc $\gamma:[0,1)\to \D$ defined by $\gamma(r):=H^{-1}\big(i(1-r)\big)$. Using Abel's equation~\eqref{EQ_Abel}, we get
$$
\lim_{r\to 1^-}\gamma(r)=\lim_{t\to +\infty}\varphi_t(z_0)=1.
$$
Hence, by Lindel\"of's Theorem (see, e.g., \cite[Theorem~9.3 on p.\,268]{Pombook75}), we have $H(1):=\angle\lim_{z\to 1}H(z)=0$.
Therefore, in order to prove (\ref{crucial}), it is enough to check that for any sequence $(a_n)\subset (0,1)$ converging to~$1$,
\begin{equation*}
\lim_{n\to\infty}\frac{1-H^{-1}(ix_n)}{|1-H^{-1}(ix_n)|}=\exp\big(i\eta\tfrac\pi2\big),\quad\text{where $x_n:=|H(a_n)|$ for all $n\in\N$.}
\end{equation*}
For such a sequence~$(a_n)$ consider the following automorphisms of~$\UD$,
$$
T_n(z):=\frac{a_n+z}{1+a_n z},\quad z\in\D,
$$
and univalent functions
\begin{equation}\label{EQ_F_n}
F_n(z):=\frac{H(T_n(z))}{x_n},\quad z\in\D.
\end{equation}
For all $n\in\Natural$, $\{it:t>0\}\subset F_n(\D)\subset \C\setminus \{-it:t\geq0\}.$
Therefore, $(F_n)$ is a normal family in~$\UD$. Since by construction $|F_n(0)|=1$ for all $n\in\Natural$, $(F_n)$ is indeed relatively compact in $\Hol(\D,\C)$. Moreover, by Hurwitz's Theorem,  any accumulation point of  $(F_n)$ is either a constant or a univalent function in~$\UD$. Let $g:\D\to\C$ be one of those accumulation points, i.e. suppose that $g$ is the limit of a some subsequence $(F_{n_k})$. Denote ${g_{k}:=F_{n_k}}$, ${k\in\Natural}$.

Since by the hypothesis, $\alpha^{-}>0$ and $\alpha^{+}>0$, there exists $\beta_1,{\beta_2>0}$ and ${\varepsilon>0}$ such that $S_0(\beta_1,\beta_2)\bigcap\,\{w:|w|<\varepsilon\}\subset H(\UD)$.  Therefore, $z=0\in\partial H(\D)$ is a well-accessible point for $H$. (For the definition of well-accessibility, we refer the reader to \cite[p.\,251]{Pombook92}.)
According to  \cite[Theorem~11.3 on p.\,251]{Pombook92}, it follows that there exist constants ${M>0}$ and ${\mu>0}$ such that for every $0\leq s\leq \rho<1$,
\begin{equation}\label{eq-pom}
|H(\rho)|=|H(\rho)-H(1)|\leq M \dist\big(H(s),\partial H(\D)\big)\left(\frac{1-\rho}{1-s}\right)^{\!\mu},
\end{equation}
where  $\dist(\cdot,\cdot)$ denotes the Euclidean distance in~$\C$; i.e., $\dist(z,W):=\inf_{w\in W}|w-z|$.

Since $0\in\partial g_k(\D)$ for all $k\in\N$, inequality~\eqref{eq-pom} with $\rho:=a_{n_k}$ and $s=s(x):=T_{n_k}(x)$ leads to
\begin{align}
\notag\dist\big(g_k(0),\partial g_k(\D\big))&\leq |g_k(0)|\le M \dist\big(g_k(x),\partial g_k(\D)\big)\left(\frac{1+a_{n_k}x}{1-x}\right)^{\!\mu} \\
\label{key1}
&\le
M\, |g_k(x)|\left(\frac{1+a_{n_k}x}{1-x}\right)^{\!\mu} \quad\text{for all $x\in[-a_{n_k},0]$.}
\end{align}
On the other hand, applying again~\eqref{eq-pom} with $\rho=\rho(x):=T_{n_k}(x)$ and ${s:=a_{n_k}}$,
we have
\begin{align}
\notag
\dist\big(g_k(x),\partial g_k(\D)\big)&\leq|g_k(x)|\leq M \dist\big(g_k(0),\partial g_k(\D)\big)\left(\frac{1-x}{1+a_{n_k}x}\right)^{\!\mu}\\
\label{key2}
&\leq M |g_k(0)|\left(\frac{1-x}{1+a_{n_k}x}\right)^{\!\mu} \quad\text{for all $x\in[0,1)$.}
\end{align}

Recall that $|g_k(0)|=1$ for all~$k\in\N$. Hence from \eqref{key1} with $x=0$, we obtain
$$
\dist\big(g_k(0),\partial g_k(\D)\big)\leq 1\leq M \dist\big(g_k(0),\partial g_k(\D)\big).
$$
Therefore,  see e.g. \cite[Theorem 3.4.9]{BCD-Book},
$$
\frac1M\leq |g_k^\prime(0)|\leq 4\quad \text{for all }~k\in\N.
$$
It follows that $g$ cannot be constant and thus it is univalent in~$\UD$. In particular, by Proposition~\ref{PR_CarathConv}, this means that the sequence of domains $D_k:=g_k(\UD)$ converges to a non-trivial kernel~$D_*$ w.r.t. $g(0)$ and that $g(\UD)=D_*$.

Denote $\Omega_r:=\frac{1}{r}\Omega$, $r>0$. On the one hand, by Theorem~\ref{kernel}, there exists $\omega\in\C$ w.r.t. which $(\Omega_r)$  converges to its kernel $\Ker{(\Omega_r)}\omega=S_0(\alpha^-,\alpha^+)$.

On the other hand, convergence of $(D_k)$ to its kernel~$D_*$ means that the sequence $(\Omega_{r_k})$, $r_k:=1/x_{n_k}$, converges to its kernel $\{w\colon-1/w\in D_*\}$ w.r.t. $\omega':=-1/g(0)$.

Using Proposition~\ref{LM_Omega-kernel}\,(C) and the definition of convergence to the kernel, we see that
$
\Ker{(\Omega_{r_k})}{\omega'}=\Ker{(\Omega_{r_k})}{\omega}=\Ker{(\Omega_r)}{\omega}=S_0(\alpha^-,\alpha^+).
$
It follows that $$g(\UD)=S_0(\alpha^+,\alpha^-).$$ Therefore, according to the Riemann Mapping Theorem,
\begin{equation}\label{EQ_g}
g(z)=i\exp\big(i\,\tfrac{\alpha^+ -\alpha^-}{2}\big)
\left(\dfrac{1-U(z)}{1+U(z)}\right)^{\!(\alpha^- +\alpha^+)/{\pi}}, \quad z\in\D,
\end{equation}
for a suitable $U\in\Aut$. We can determine~$U$ using~\eqref{key1} and~\eqref{key2}. Indeed, passing in these inequalities to the limit as ${k\to+\infty}$ and taking into account that $|g_k(0)|=1$ for all~$k\in\Natural$ and that $\lim\limits_{n\to+\infty}a_{n}=1$, we get
$$
|g(x)|\ge \frac1M\left(\frac{1-x}{1+x}\right)^{\!\mu}\text{~for all $x\in(-1,0]$}~
\text{~and~}~ |g(x)|\le M\left(\frac{1-x}{1+x}\right)^{\!\mu}\text{~for all $x\in[0,1)$.}
$$
It follows that ${g(x)\to\infty}$ as ${x\to-1^+}$ and ${g(x)\to0}$ as ${x\to1^-}$.
Taking into account that $|g(0)|=1$ and using~\eqref{EQ_g}, we therefore conclude that $U=\id_\UD$.

We have proved that every converging subsequence of~$(F_n)$ has the same limit. Recalling that $(F_n)$ is a normal family in~$\UD$, we may conclude that $(F_n)$ converges locally uniformly in $\UD$ to
\begin{equation}\label{EQ_F}
F(z):=i\exp\big(i\,\tfrac{\alpha^+ -\alpha^-}{2}\big)
\left(\dfrac{1-z}{1+z}\right)^{\!(\alpha^- +\alpha^+)/{\pi}}, \quad z\in\D.
\end{equation}

Note that $i\in F_n(\UD)$ for all~$n\in\Natural$ and that $i\in F(\UD)$. Hence by Proposition~\ref{PR_CarathConv},
\begin{equation*}
z_n:=F^{-1}_{n}(i)\to z_0:=F^{-1}(i)\in\UD\quad \text{as~$n\to+\infty$}.
\end{equation*}
Furthermore, by~\eqref{EQ_F_n} with $z:=z_n$, for all $n\in\Natural$, we have
$$
{1-H^{-1}(ix_{n})}={1-T_n(z_n)}={(1-a_n)(1-z_n)/(1+a_nz_n)}.
$$
Therefore,
$$
\frac{1-H^{-1}(ix_n)}{|1-H^{-1}(ix_n)|}=\frac{1-z_n}{|1-z_n|}\frac{|1+a_nz_n|}{1+a_nz_n}\to
\frac{1-z_0}{1+z_0}\frac{|1+z_0|}{|1-z_0|}~\text{~as $n\to+\infty$}.
$$
Finally, according to~\eqref{EQ_F}, we have
$$
\frac{1-z_0}{1+z_0}=\Big(\exp\big(i\,\tfrac{\alpha^- - \alpha^+}{2}\big)\Big)^{\pi/(\alpha^- +\alpha^+)}=\exp\big(i\eta\tfrac\pi2\big).
$$
This completes the proof.\qed
\end{proof}
\end{theopargself}

Now, we are going to apply the above results to  domains starlike at infinity whose boundary is contained in a ``neighbourhood'' of the boundary of a sector $S_p(\beta_1,\beta_2)$.

\begin{corollary}\label{CR1}
Let $(\varphi_t)$ be a non-elliptic semigroup in $\D$ with DW-point $\tau\in\partial\D$, Koenigs function $h$, and planar domain $\Omega:=h(\D)$. Let $\rho:\Complex\to[0,+\infty)$ be a continuous function such that ${\rho(w)/|w|\to0}$ as ${w\to\infty}$. Fix some $p\in\C$ and $\beta_1,\beta_2\in[0,\pi]$ with $\beta_1+\beta_2>0$ and suppose that
\begin{equation}\label{EQ_boundary-of-Omega}
\dist\big(w,\partial S_p(\beta_1,\beta_2)\big)\le\rho(w) \quad\text{for any $~w\in\partial\Omega$}.
\end{equation}
If $\beta_1\beta_2=0$, we additionally require that for all $R>0$ large enough,
\begin{align}\label{EQ_boundary-additionalA}
&p+iR\,\,\exp\big(i(\beta_1-\beta_2)/2\big)\in\Omega,\\
\intertext{and if $\beta_1=\pi\neq\beta_2$ or $\beta_1\neq\pi=\beta_2$,  then we additionally require that}
\label{EQ_boundary-additionalB}
&p-iR\,\exp\big(i(\beta_1-\beta_2)/2\big)\in\C\setminus\Omega\quad\text{for all $R>0$ large enough}.
\end{align}

Then for all $z\in\D$,
$$
\mathrm{Slope}[t\mapsto\varphi_t(z),\tau]=\big\{\eta\tfrac{\pi}{2}\big\},\quad \eta:=\frac{\beta_1-\beta_2}{\beta_1+\beta_2}.
$$
\end{corollary}
\begin{remark}
Note that additional conditions~\eqref{EQ_boundary-additionalA} and~\eqref{EQ_boundary-additionalB} in Corollary~\ref{CR1} cannot be omitted. For example, if $\Omega$ satisfies condition~\eqref{EQ_Betsakos-cond} below, then it satisfies also~\eqref{EQ_boundary-of-Omega}~--- for a suitable $p\in\C$ and a constant function~$\rho$~--- whenever $\pi\in\{\beta_1,\beta_2\}$. Condition~\eqref{EQ_boundary-additionalB} excludes all the cases, except for $\beta_1=\beta_2=\pi$. Another similar example is provided by any hyperbolic one-parameter semigroup, for which the conclusion of Corollary~\ref{CR1} is not valid (see Remark~\ref{RM_slopes}). Since the planar domain of a hyperbolic semigroup is contained in some vertical strip, condition~\eqref{EQ_boundary-of-Omega} would be satisfied for such a semigroup both with $(\beta_1,\beta_2):=(0,\pi)$ and with $(\beta_1,\beta_2):=(\pi,0)$. At the same time, conditions~\eqref{EQ_boundary-additionalA} and~\eqref{EQ_boundary-additionalB} fail in this case.
\end{remark}

In the special case $\rho=\const$ and $\beta_1=\beta_2=\pi$, we recover a result of Betsakos~\cite{Bet16a}.
\begin{corollary}[\protect{\cite[Theorem 2]{Bet16a}}, see also \protect{\cite[Corollary~5.1\,(3)]{BCDG19}}]\label{CR2}\\
Let $(\varphi_t)$ be a non-elliptic semigroup in $\D$  with DW-point $\tau\in\partial\D$ and Koenigs function~$h$ and let $\Omega:=h(\D)$. If there exist positive numbers $a_1$, $a_2$, and $b$ such that
\begin{equation}\label{EQ_Betsakos-cond}
\partial\Omega\subset \{x+iy: a_1<x<a_2,\ y<b\},
\end{equation}
then for all $z\in\UD$,
$$
\mathrm{Slope}[t\mapsto\varphi_t(z),\tau]=\left\{0\right\}.
$$	
\end{corollary}
Assuming now that $(\beta_1,\beta_2)\neq(\pi,\pi)$, for $\rho=\const$ we obtain the following statement.
\begin{corollary}\label{CR3}
Let $(\varphi_t)$ be a non-elliptic semigroup in $\D$  with DW-point $\tau\in\partial\D$, Koenigs function~$h$, and planar domain $\Omega:=h(\D)$. Fix arbitrary $\beta_1,\beta_2\in[0,\pi]$ with $0<\beta_1+\beta_2<2\pi$. If for some ${p\in\C}$ and some $q\in S_p(\beta_1,\beta_2)$,
\begin{equation*}\label{EQ_BCDG-cond}
S_{q}(\beta_1,\beta_2)\,\subset\,\Omega\,\subset\, S_p(\beta_1,\beta_2),
\end{equation*}
Then for all $z\in\D$,
$$
\mathrm{Slope}[t\mapsto\varphi_t(z),\tau]=\big\{\eta\tfrac{\pi}{2}\big\},\quad \eta:=\frac{\beta_1-\beta_2}{\beta_1+\beta_2}.
$$
\end{corollary}

Setting $\beta_1=\beta_2$ in the above corollary, we immediately obtain the statements (1) and~(2) of \cite[Corollary~5.1]{BCDG19}.

Since Corollaries~\ref{CR2} and~\ref{CR3} follow directly from Corollary~\ref{CR1}, we only need to prove the latter one.
Two examples making use of Corollary~\ref{CR1}  with a non-constant function~$\rho$ can be found at the beginning of Section~\ref{S_examples}.

\begin{theopargself}
\begin{proof}[\textit{of Corollary~\ref{CR1}}~]
Clearly, without loss of generality we may assume that $p=0$.

Fix some $\theta\in(0,1)$. Denote $\zeta:=i\exp\big(i(\beta_1-\beta_2)/2\big)$. The ray $\{R\zeta\colon R>0\}$ is the internal bisector of $S_0(\beta_1,\beta_2)$. Hence
$$
 A_R:=\big\{te^{i\psi}\zeta\colon t>R,~|\psi|<\theta(\beta_1+\beta_2)/2\big\}~\subset~ S_0(\beta_1,\beta_2)
$$
for any~${R>0}$. Moreover, there exists~$\varepsilon>0$ such that
\begin{equation}\label{EQ_boundary-sectors}
\hphantom{\dist}\dist\big(w,\partial S_0(\beta_1,\beta_2)\big)>\varepsilon |w| \quad\text{for any $~R>0~$ and all $~w\in A_R$}.
\end{equation}
Taking into account that $\dist\big(0,A_R\big)=R$, we see that there exists ${R_0>0}$ such that $\rho(w)\le \varepsilon |w|$ for all $w\in A_{R_0}$. Thanks to \eqref{EQ_boundary-of-Omega} and~\eqref{EQ_boundary-sectors}, it follows that $A_R\bigcap\,\partial\Omega=\emptyset$ and hence, either ${A_{R_0}\subset\Omega}$ or ${A_{R_0}\subset\Complex\setminus\Omega}$.

Consider the following cases.\medskip

\noindent\textit{Case 1:} $\beta_1=\beta_2=\pi$. In this case, for any $w\in\Complex$, the ray $\{w+it\colon t\ge0\}$ intersects~$A_{R_0}$. Hence we may conclude that ${A_{R_0}\subset\Omega}$. It follows that $0\in\mathsf{NP}(\Omega)$ and
$$
 \theta \pi\le \alpha_{\Omega,0}^-(t),\alpha_{\Omega,0}^+(t)\le\pi\quad\text{for all $t>R_0$}.
$$
Since $\theta$ can be chosen as close to~$1$ as we wish, this means that $\alpha_{\Omega,0}^-(t),\alpha_{\Omega,0}^+(t)\to\pi$ as ${t\to+\infty}$ and it remains to refer to Theorem~\ref{main}.

\medskip

From now on we will suppose that $(\beta_1,\beta_2)\neq(\pi,\pi)$. Arguing as above, we see that for any $\theta\in(0,1)$, there exists $R_1>0$ such that $B_{R_1}\bigcap\,\partial\Omega=\emptyset$, where
$$
B_{R}:=\big\{-te^{i\psi}\zeta\colon t>R,~|\psi|<\theta(2\pi-\beta_1-\beta_2)/2\big\}.
$$

\noindent\textit{Case 2:} $\beta_1,\beta_2\in(0,\pi)$. If $\theta$ is sufficiently close to~$1$, then for any $w\in\C$,
$$
  \textstyle\{w+it\colon t\ge0\}\,\bigcap\,A_{R_0}\neq\emptyset\quad\text{and}\quad \{w+it\colon t\le0\}\,\bigcap\,B_{R_1}\neq\emptyset.
$$
It follows that $A_{R_0}\subset\Omega$ and $B_{R_1}\subset\C\setminus\Omega$. Therefore, $0\in\mathsf{NP}(\Omega)$ and for every \\$t>\max\{R_0,R_1\}$ we have
\begin{gather*}
\beta_1-(1-\theta)\frac{\beta_1+\beta_2}2~\le~\alpha_{\Omega,0}^-(t) ~\le~\beta_1+(1-\theta)\frac{2\pi-\beta_1-\beta_2}2\\
\mathllap{\text{and}\qquad}
\beta_2-(1-\theta)\frac{\beta_1+\beta_2}2~\le~\alpha_{\Omega,0}^+(t) ~\le~\beta_2+(1-\theta)\frac{2\pi-\beta_1-\beta_2}2.
\end{gather*}
Again, in this case, the conclusion of the corollary follows from Theorem~\ref{main}.

\medskip\noindent\textit{Case 3:} $0<\beta_1<\pi$, $\beta_2=\pi$. As in the previous case, we see that $A_{R_0}\subset\Omega$. Thanks to~\eqref{EQ_boundary-additionalB}, we also have $B_{R_1}\subset\C\setminus\Omega$. The rest of the proof is the same as in Case~2.

\medskip\noindent\textit{Case 4:}
$\beta_1=0$, $0<\beta_2<\pi$. As in Case~2, we see that $B_{R_1}\subset\C\setminus\Omega$. Moreover, condition~\eqref{EQ_boundary-additionalA} allows us to conclude that $A_{R_0}\subset\Omega$.
Since $\Omega$ is starlike at infinity,  the latter inclusion implies that there exists $q\in\Omega$ such that $S_q(0,\theta\beta_2)\subset\Omega$. It follows that
$$
\liminf_{t\to+\infty}\alpha_{\Omega,q}(t)\ge\beta_2>0.
$$
Moreover, since $\theta$ can be chosen arbitrarily close to~$1$, the inclusion $B_{R_1}\subset\C\setminus\Omega$ implies that $\alpha_{\Omega,q}^-(t)\to 0$ as $t\to+\infty$. Therefore, by Proposition~\ref{PR_second}\,(C), $\mathrm{Slope}[t\mapsto \varphi_t(z),\tau]=\{-\pi/2\}$.

\medskip\noindent\textit{Case 5:} $\beta_1=0$, $\beta_2=\pi$. Conditions~\eqref{EQ_boundary-additionalA} and~\eqref{EQ_boundary-additionalB} allows us to conclude that $A_{R_0}\subset\Omega$ and $B_{R_1}\subset\C\setminus\Omega$. As in the previous case, using the fact that $\Omega$ is starlike at infinity, we see that there exists $q\in\Omega$ such that $S_q(0,\theta\pi)\subset\Omega$. The rest of the proof is literally the same as in Case~4.

\medskip We omit the remaining three cases: ${\beta_1=\pi}$ and ${0<\beta_1<\pi}$; ${0<\beta_1<\pi}$ and ${\beta_2=0}$; ${\beta_1=\pi}$ and ${\beta_2=0}$, because they are analogous to Cases 3, 4, and 5, respectively.
\qed
\end{proof}
\end{theopargself}

\label{seq}Denote by $\mathfrak S$ the set of all sequences $(t_n)\subset(0,+\infty)$ tending to~${+\infty}$ and such that $${\sup\limits_{n\in\Natural}|t_{n+1}-t_n|<+\infty}.$$
\begin{definition} Let $\Omega$ be a domain starlike at infinity, $p\in\C$ and $\alpha,\beta\in[0,\pi]$ with $\alpha+\beta>0$. We say that $\Omega$ \textsl{meets $S_p(\alpha,\beta)$ on the left} (resp. \textsl{on the right}) \textsl{at uniform times} if there exists a sequence $(t_n)\in\mathfrak{S}$  such that
\begin{equation}
\{p+it_ne^{i\alpha}:n\in\N\}\subset \partial\Omega\quad \text{(resp. $\{p+it_ne^{-i\beta}:n\in\N\}\subset \partial\Omega$)}.
\end{equation}
\end{definition}

\begin{corollary}\label{CR_meet} Let $(\varphi_t)$ be a non-elliptic semigroup in $\D$ with DW-point $\tau\in\partial\D$, Koenigs function $h$, and planar $\Omega:=h(\D)$. Let $p\in\C$. The following statements hold.
\begin{enumerate}
\item[(A)] Assume there exist $\beta_1,\beta_2\in(0,\pi)$ such that
$S_p(\beta_1,\beta_2)\subset \Omega$
and $\Omega$ meets $S_p(\beta_1,\beta_2)$ on the right and on the left at uniform times. Then
$$
 \mathrm{Slope}[t\mapsto\varphi_t(z),\tau]=\big\{\eta\tfrac{\pi}{2}\big\},~ \eta:=\frac{\beta_1-\beta_2}{\beta_1+\beta_2},~\text{~for all ${z\in\UD}$.}
$$

\item[(B)] Assume there exist $\beta\in(0,\pi)$ such that
$
S_p(\pi,\beta)\subset \Omega,
$
and $\Omega$ meets $S_p(\pi,\beta)$ on the right at uniform times. Then for all ${z\in\UD}$,
$$
 \mathrm{Slope}[t\mapsto\varphi_t(z),\tau]=\big\{\eta\tfrac{\pi}{2}\big\},\quad \eta:=\frac{\pi-\beta}{\pi+\beta}.
$$

\item[(C)] Assume there exist $\beta\in(0,\pi)$ such that
$
S_p(\beta,\pi)\subset \Omega,
$
and $\Omega$ meets $S_p(\pi,\beta)$ on the left at uniform times. Then for all ${z\in\UD}$,
$$
 \mathrm{Slope}[t\mapsto\varphi_t(z),\tau]=\big\{\eta\tfrac{\pi}{2}\big\},\quad \eta:=\frac{\beta-\pi}{\beta+\pi}.
$$
\end{enumerate}	
\end{corollary}
\begin{proof}
The hypothesis of~(A) implies that $p\in\mathsf{NP}(\Omega)$ and that
$$
 \alpha^-_{p,\Omega}(t'_n)=\beta_1~\text{ and }~ \alpha^+_{p,\Omega}(t''_n)=\beta_2,\quad n\in\Natural,
$$
for some sequences $(t_n'), (t_n'')\in\mathfrak{S}$. Note that $t'_{n+1}/t'_n,\,t''_{n+1}/t''_n\to1$ as ${n\to+\infty}$. Therefore, by Theorem~\ref{kernel},
${\alpha^-_{\Omega,p}(t)\to\beta_1}$ and ${\alpha^+_{\Omega,p}(t)\to\beta_2}$ as ${t\to+\infty}$. Thus, the desired conclusion holds by Theorem~\ref{main}.
\medskip

\noindent\textit{Proof of (B)}  Since $S_p(\pi,\beta)\subset\Omega$, we have $p\in\mathsf{NP}(\Omega)$ and $\alpha_{\Omega,p}^-(t)=\pi$ for all ${t>0}$. Moreover, since $\Omega$ meets $S_p(\pi,\beta)$ on the right at uniform times, $\alpha_{\Omega,p}^-(t_n)=\beta$ for a suitable sequence $(t_n)\in\mathfrak{S}$. Therefore, as above, the desired conclusion follows from Theorems~\ref{kernel} and~\ref{main}.
\medskip

\noindent\textit{Proof of (C)} is omitted because it is similar to that of assertion~(B).\qed
\end{proof}

\section{Examples}\label{S_examples}
We start this section with a few simple examples illustrating Corollary~\ref{CR1}. Recall that any domain $\Omega$ starlike at infinity and different from~$\C$ defines a non-elliptic one-parameter semigroup (Remark~\ref{RM_Koenings-image}). Moreover, this semigroup is parabolic and of zero hyperbolic step if and only if $\mathsf{NP}(\Omega)=\C$.

\begin{example}
Let $f:\Real\to\Real$ be a continuous function such that the limits
$$
 \kappa_1:=\lim_{x\to-\infty}\frac{f(x)}{|x|},\quad \kappa_2:=\lim_{x\to+\infty}\frac{f(x)}{x}
$$
exist and they are finite.
Then for the parabolic one-parameter semigroup $(\varphi_t)$ with zero hyperbolic step whose planar domain is
$
\Omega:=\{x+iy\colon y>  f(x)\}
$
we have $$
\mathrm{Slope}[t\mapsto\varphi_t(z),\tau]=\big\{\eta\tfrac\pi2\big\},\quad
\eta:=\frac{\arctan k_1-\arctan k_2}{\arctan k_1+\arctan k_2},
$$
for any~$z\in\UD$. Indeed, the hypothesis of Corollary~\ref{CR1} is satisfied in this case with $p:=0$, $\beta_j:=\arctan \kappa_j$, $j=1,2$, and
$$
\rho(x+iy):=\begin{cases}
|f(x)-\kappa_1x| &\text{if~$x<0$},\\[1ex]
|f(x)-\kappa_2x| &\text{if~$x\ge0$}.
\end{cases}
$$
\end{example}

\begin{example}
Corollary~\ref{CR1} can be applied for a one-parameter semigroup with planar domain $\Omega:=\{x+iy\colon y>x^3\}$ if we set $\rho(w):=|w|^{1/3}$ and $p:=0$. This example illustrates the role of conditions \eqref{EQ_boundary-additionalA} and~\eqref{EQ_boundary-additionalB}: in this case, \eqref{EQ_boundary-of-Omega} is satisfied both with $(\beta_1,\beta_2):=(\pi,0)$ and with $(\beta_1,\beta_2):=(0,\pi)$; however, conditions~\eqref{EQ_boundary-additionalA} and~\eqref{EQ_boundary-additionalB} exclude the latter possibility.
\end{example}

\begin{example}
Let $(t'_n),(t''_n)\in\mathfrak{S}$, see page~\pageref{seq}. Fix some $\beta_1,\beta_2\in(0,\pi)$ and let
$$
E_1:=\bigcup_{n\in\Natural}\{it'_ne^{i\beta_1}+iy\colon y\le0\},\quad E_2:=\bigcup_{n\in\Natural}\{it''_ne^{-i\beta_2}+iy\colon y\le0\}.
$$
Then by Corollary~\ref{CR_meet}, $\Omega_1:=\C\setminus E_1$, $\Omega_2:=\C\setminus E_2$, and $\Omega_3:=\C\setminus\big(E_1\bigcup E_2\big)$ are the planar domains of parabolic one-parameter semigroups $(\varphi^k_t)$, $k=1,2,3$, respectively, with
$$
 \mathrm{Slope}[t\mapsto\varphi^k_t(z),\tau]=\big\{\eta_k\tfrac\pi2\big\}\quad\text{for all $z\in\UD$},
$$
where
$$
\eta_1:=\frac{\beta_1-\pi}{\beta_1+\pi},\quad\eta_2:=\frac{\pi-\beta_2}{\pi+\beta_2},\quad \eta_3:=\frac{\beta_1-\beta_2}{\beta_1+\beta_2}.
$$
\end{example}

The next two examples illustrate the difference between the distance functions $\delta^\pm_{\Omega,p}$ and the angular extents $\alpha^\pm_{\Omega,p}$ as characteristics of the geometry of a domain  starlike at infinity $\Omega\subsetneqq\C$ near the point~$\infty$.

\begin{example}\label{EX_Santiago}
There exists a domain $\Omega$ starlike at infinity and a point $p\in \mathsf{NP}(\Omega) $ such that the functions $t\mapsto \delta^{+}_{\Omega,p} (t)$ and $t\mapsto t\alpha^{+}_{\Omega,p} (t)$ are not asymptotically equivalent as~${t\to+\infty}$.
\end{example}

\begin{proof} For $n\in\N$ we set $t_{n}:=n!$,  ${\gamma:=\sqrt{2}/2}$, $y_{n}:=\sqrt{t_{n+1}^{2}-t_{n}^2\gamma^{2}}$,
\begin{align*}
\Gamma_{0}&:=\{ w\in \C:\ \Re w= 0,\ \Im w\leq 0\},\\
\text{and}\quad\Gamma_{n}&:=\{ w\in \C:\ \Re w= t_{n} \gamma,\ \Im w\leq y_{n}-1\}.
\end{align*}\
Consider the domain
\begin{equation}\label{Eq:domain}
\Omega:=\C~\big\backslash \bigcup_{n=0}^{\infty}\Gamma_{n}.
\end{equation}
Notice that for each $n\in\N$,
the point $t_n\gamma+iy_n$ lies on the semicircle $$C_{n+1}:=\{w\colon |w|=t_{n+1},\,\Im w\ge0\}.$$ It follows that the slits $\Gamma_{m}$ with ${m\leq n}$ do not intersect~$C_{n+1}$. Moreover, $t_k\gamma> t_{n+1}$ for any ${k>n+1}$. Therefore, the slits $\Gamma_k$ with ${k> n+1}$ do not intersect~$C_{n+1}$ either.
On the other hand,
\begin{equation}\label{EQ_yt}
 y_{n+1}-1\,=\,t_{n+1}\sqrt{(n+2)^2-\gamma^2}\,-\,1\,>\,t_{n+1}\quad\text{for all~$~n\in\N$}.
\end{equation}
It follows that $C_{n+1}$ intersect $\Gamma_{n+1}$ at the point $t_{n+1}(\gamma+i\sqrt{1-\gamma^2})=t_{n+1}e^{i\pi/4}$.
Thus
\begin{equation}\label{EQ_alpha}
 \alpha^{+}_{\Omega,0} (t_{n+1})=\pi/4\quad\text{ for all $~n\in\N$.}
\end{equation}

On the other hand, by the very definition,
\begin{equation*}
\delta^{+}_{\Omega,0} (t_{n+1})\le |it_{n+1}-t_n\gamma-i(y_n-1)|.
\end{equation*}
Using the triangle inequality, we obtain
$$
|it_{n+1}-t_{n}\gamma-i(y_{n}-1)|~\le~|it_{n+1}-t_{n}\gamma-iy_{n}|+1=\sqrt{2t_{n+1}(t_{n+1}-y_n)}+1
$$
for any $n\in\Natural$. Furthermore,
$$
 \theta_n:=\frac{t_{n+1}-y_n}{t_{n+1}}\cdot 4(n+1)^2=\left(1-\sqrt{1-\frac1{2(n+1)^2}}\,\right)\cdot4(n+1)^2~\to~1
$$
as $\,n\to+\infty$.  Hence,
$$
\delta^{+}_{\Omega,0} (t_{n+1})~\le~ |it_{n+1}-t_{n}\gamma-i(y_{n}-1)|~\le~\sqrt{\theta_n}
\frac{ t_{n+1}}{\sqrt2(n+1)}+1.
$$

Thus, taking into account~\eqref{EQ_alpha}, we have
$$
\lim_{n\to +\infty}\frac{\delta^{+}_{\Omega,0}(t_{n+1})}{\,t_{n+1}\alpha^{+}_{\Omega,0}(t_{n+1})}=0.
$$
In particular, $\delta^{+}_{\Omega,0}(t)$ and $t\mapsto t\alpha^{+}_{\Omega,0}(t)$ are not asymptotically equivalent at~$+\infty$.\qed
\end{proof}

The next example shows that it is not possible to get a result similar to Theorem~\ref{maintheoremBCDGZ} using the functions $\alpha^\pm_{\Omega,p}$ instead of $\delta_{\Omega,p}^\pm $.
\begin{example}\label{EX_Manolo} There exists a parabolic semigroup $(\varphi_t)$ in $\D$ of zero hyperbolic step with the associated planar domain~$\Omega$ and a sequence  $(t_n)\subset(0,+\infty)$ tending to $+\infty$ such that $(\varphi_{t_{n}}(z))$ converges to the DW-point of the semigroup \textit{non-tangentially}, but $\alpha^-_{\Omega,0}(t_{n})$ and $\alpha^+_{\Omega,0}(t_{n})$ are \textit{not} asymptotically equivalent.
\end{example}

\begin{proof} For $n\in\Natural$, we denote $t_n:=n!$, $\alpha_{n}:=\arcsin(1/n)\in (0,\pi/2]$, and $$y_{n}:=\sqrt{t_{n+1}^{2}-(t_{n}\sin\alpha_n)^2}.$$
Furthermore, for each $n\in \N$, let
\begin{equation}\label{Eq:slits}
\begin{split}
&\Gamma_{n}:=\{ w\in \C:\ \Re w= t_{n}\sin\alpha_n,\ \Im w\leq y_{n}-1\},\\
&\Lambda_{n}:=\{ w\in \C:\ \Re w=-t_{n}\sin\alpha_n,\ \Im w\leq y_{n}\},\\
&\Gamma :=\{ w\in \C:\ \Re w= 0,\ \Im w\leq 0\}.
\end{split}
\end{equation}
Consider the domain
\begin{equation}\label{Eq:domain}
\Omega:=\C~\Big\backslash \Big(\Gamma\,\cup\,\bigcup\limits_{n=2}^{+\infty}(\Gamma_{n}\cup\Lambda_{n})\Big)
\end{equation}
sketched in Figure~\ref{Fig1}. Clearly $\Omega$ is starlike at infinity. Fix any conformal map $h$ of~$\UD$ onto~$\Omega$ and consider the semigroup $(\varphi_t)$ defined by $\varphi_{t}:=h^{-1}\circ(h+it)$ for all ${t\ge0}$. To simplify the notation we write $\delta^\pm(t):=\delta_{\Omega,0}^\pm (t)$ and $\alpha^\pm(t):= \alpha^\pm_{\Omega,0} (t)$ for all $t>0$.

 \begin{figure}



\begin{center}
\begin{tikzpicture}[scale=0.6]


\draw [dashed]  (0,-1)--(0,25);
\draw [dashed] (-6,0)--(14,0);


\draw (0,24) node[scale=1] {$\bullet$};
\draw (0.5,24.5) node[scale=1] {$it_{4}$};

\draw (0,6) node[scale=1] {$\bullet$};
\draw (0.5,6.5) node[scale=1] {$it_{3}$};

\draw (0,2) node[scale=1] {$\bullet$};
\draw (0.5,2.5) node[scale=1] {$it_{2}$};


\draw  [thick]  (0,-1)--(0,0);
\draw (0,-1.5) node[scale=1] {$\Gamma$};

\draw  [thick]  (12,-1)--(12,24);
\draw  [thick, dashed]  (12,24)--(12,25);
\draw (12,-1.5) node[scale=1] {$\Gamma_{4}$};

\draw  [thick]  (3.46,-1)--(3.46,22.75);
\draw (3.46,-1.5) node[scale=1] {$\Gamma_{3}$};

\draw  [thick]  (1.4,-1)--(1.4,4.83);
\draw (1.4,-1.5) node[scale=1] {$\Gamma_{2}$};



\draw [thick] (-3.46,-1)--(-3.46,23.75);
\draw (-3.46,-1.5) node[scale=1] {$\Lambda_{3}$};
\draw  [thick]  (-1.4,-1)--(-1.4,5.83);
\draw (-1.4,-1.5) node[scale=1] {$\Lambda_{2}$};


\draw [dashed] [domain=0:12] plot (\x,{0.025+sqrt(576-(\x)^2)});
\draw [dashed] [domain=-3.2:-0.1] plot (\x,{0.025+sqrt(576-(\x)^2)});

\draw [dashed] [domain=-1.4:3.46] plot (\x,{0.015+sqrt(36-(\x)^2)});


\draw (6,24) node[scale=1] {$\alpha^{+}_{\Omega,0}(t_{4})$};


\draw (-2,24.5) node[scale=1] {$\alpha^{-}_{\Omega,0}(t_{4})$};


\draw[{Stealth[scale=1.25,
          length=6,
          width=4]}-{Stealth[scale=1.25,
          length=6,
          width=4]}] (0.05,23.97) -- (3.46,22.75);

\draw (1.8,22.2) node[scale=1] {$\delta^{+}_{\Omega,0}(t_{4})$};

\draw[{Stealth[scale=1.25,
          length=6,
          width=4]}-{Stealth[scale=1.25,
          length=6,
          width=4]}] (-0.05,23.99) -- (-3.48,23.75);

\draw (-1.9,23) node[scale=1] {$\delta^{-}_{\Omega,0}(t_{4})$};

%
%
%
%
%
%
%
%
%
%
%
%
%
%
%
%
%
%
%

\end{tikzpicture}
\end{center} 
  \caption{}\label{Fig1}
\end{figure}

Let us show that $(\varphi_{t_{n}}(z))$ converges non-tangentially to the DW-point of the semigroup~$(\varphi_t)$. Bearing in mind Theorem \ref{maintheoremBCDGZ}, we have to prove that $\delta^+(t_{n})$ behaves asymptotically like $\delta^-(t_{n})$ as ${n\to+\infty}$.

Fix for a while some $n\ge2$ and $\kappa\in\{0,1\}$. For $m\in\{1,\ldots,n\}$ denote $c_m:=it_{n+1}-w_m$, where $w_m:={(-1)^{1-\kappa}t_m\sin\alpha_m+iy_m-i\kappa}$ is the tip of the slit~$\Lambda_m$ if~${\kappa=0}$ or that of~$\Gamma_m$ if~${\kappa=1}$.
Using the inequality $\sqrt{x}>x$ valid for all~$x\in(0,1)$ and taking into account that ${y_m<t_{m+1}}$, we see that
\begin{equation}\label{EQ_Delta-y}
y_{m+1}-y_m \, > y_{m+1}-t_{m+1}\, > (m+1)\,(m+1)!\,-\,(m+1)!\,=\,m(m+1)!
\end{equation}
for all $m\in\N$. Hence if $1\le m<n$, then
\begin{align*}
(\Im c_m)^2-(\Im c_{m+1})^2 \,&=(y_{m+1}-y_m)\big(2t_{n+1}-y_{m}-y_{m+1}+2\kappa\big)\\&>\, m(m+1)!\cdot\big(2t_{n+1}-t_{m+1}-t_{m+2}\big)\\&\ge\, m(m+1)!\cdot\big((n+1)!-n!\big)\,=\,mn(m+1)!\,n!,
\end{align*}
from which it follows that
\begin{multline*}
|c_m|^2\,-\,|c_{m+1}|^2~>~mn(m+1)!\,n!\,-\,
\big((m!)^2-((m-1)!)^2\big)\\
=~mn(m+1)!\,n!\,-\,(m^2-1)((m-1)!)^2~\ge~0
\end{multline*}
whenever $1\le m<n$. Therefore, for all such~$m$, we have $ |c_m|>|c_n|$ and
\begin{align*}
|c_n|^2~&=~(t_{n+1}-y_n+\kappa)^2+(t_n\sin\alpha_n)^2\\[.3ex]&=~
\Big(t_{n+1}-t_{n+1}\sqrt{1-\big(\tfrac{t_{n}}{t_{n+1}}\sin\alpha_n\big)^2}+\kappa\Big)^2+
(t_n\sin\alpha_n)^2\\[.3ex]
&<~\big((n-1)!\big)^2\Big(\tfrac1{n(n+1)}+\tfrac1{(n-1)!}\Big)^2+\big((n-1)!\big)^2\\[.3ex]
&<~3\big((n-1)!\big)^2\, <~(t_{n+1}\sin\alpha_{n+1})^2\,<\,t^2_{n+1},
\end{align*}
where we again used the fact that $\sqrt{x}>x$ if~$0<x<1$.
With $\kappa:=1$ it follows that
\begin{align}
\notag
\delta^+(t_{n+1})^{2}~&=~|it_{n+1}-(t_{n}\sin\alpha_{n}+iy_{n}-i)|^{2}=(t_{n+1}-y_{n}+1)^{2}+(t_n\sin\alpha_n)^2\\
\label{EQ_delta-plus}
&=~2t_{n+1}^{2}-2t_{n+1}y_{n}+2t_{n+1}-2y_{n}+1~\text{~for any~$~n\ge2$,}\\
\intertext{and with $\kappa:=0$ we get}
\notag
\delta^-(t_{n+1})^{2}~&=~|it_{n+1}-(-t_{n}\sin\alpha_{n}+iy_{n})|^{2}=(t_{n+1}-y_{n})^{2}+(t_n\sin\alpha_n)^2\\
\label{EQ_delta-minus}
&=~2t_{n+1}^{2}-2t_{n+1}y_{n}~\text{~for any~$~n\ge2$.}
\end{align}
Note that
\begin{equation}\label{EQ_tn-yn}
t_{n+1}-y_n\ge t_{n+1}-t_{n+1}\Big(1-\frac{(t_{n}\sin\alpha_n)^2}{2t_{n+1}^{2}}\Big)=\frac{(n-1)!}{2n(n+1)},
\end{equation}
where we used the inequality $\sqrt{1+x}\le1+x/2$ for all $x\ge-1$.
Combining \eqref{EQ_delta-plus}, \eqref{EQ_delta-minus}, and \eqref{EQ_tn-yn}, we see that
\begin{equation*}
\frac{\delta^+(t_{n+1})^{2}}{\delta^-(t_{n+1})^{2}}=1+\frac{1}{t_{n+1}}+\frac{1}{2t_{n+1}(t_{n+1}-y_{n})}\to1 \quad\text{as~$~n\to+\infty$}.
\end{equation*}

Therefore, on the one hand by Theorem \ref{maintheoremBCDGZ}, the sequence $(\varphi_{t_n}(z))$ converges to the DW-point non-tangentially.
On the other hand, by~\eqref{EQ_Delta-y}, $y_{1}<y_2<\ldots <y_n$ and $t_{n+1}<y_{n+1}-1$. Hence, by the construction,
$$
 t_{n+1}\sin\,\alpha^{+}(t_{n+1})=t_{n+1}\sin\,\alpha_{n+1}\quad\text{and}\quad t_{n+1}\sin \alpha^{-}(t_{n+1})= t_{n}\sin \, \alpha_n
$$
for all $~n\ge2$. Therefore,  $$\frac{\sin \alpha^{-}(t_{n+1})}{\sin \alpha^{+}(t_{n+1})}= \frac{t_{n}\,\sin\alpha_{n}}{t_{n+1}\,\sin\alpha_{n+1}}=\frac1{n}.$$
In particular, $\alpha^{-}(t_{n+1})/\alpha^{+}(t_{n+1})\to0$ as $n\to+\infty$.\qed
\end{proof}


\begin{thebibliography}{}
\bibitem{Abate} Abate, M.: \textsl{Iteration Theory of Holomorphic Maps on Taut
Manifolds.} Mediterranean Press, Rende, Cosenza, 1989.

\bibitem{AroBra16} Arosio, L., Bracci, F.: Canonical models for holomorphic iteration. Trans. Amer. Math. Soc. \textbf{368} (2016), no.\,5, 3305--3339.

\bibitem{Bet16a} Betsakos, D.: On the asymptotic behavior of the trajectories of semigroups of holomorphic functions. J. Geom. Anal.,  \textbf{26} (2016), 557--569.

\bibitem{BerPor78} Berkson, E., Porta, H.: Semigroups of
holomorphic functions and composition operators. Michigan
Math. J. \textbf{25} (1978), 101--115.

\bibitem{BCDG19} Bracci, F., Contreras, M.D., D\'{\i}az-Madrigal, S., Gaussier, H.: A characterization of orthogonal convergence in simply connected domains, J. Geom. Anal. {\bf 29} (2019), no.~4, 3160--3175. 

\bibitem{BCDG20} Bracci, F., Contreras, M.D., D\'{\i}az-Madrigal, S., Gaussier, H.: Non-tangential limits and the slope of trajectories of holomorphic semigroups of the unit disc. Trans. Amer. Math. Soc. {\bf 373} (2020), no.\,2, 939--969. 

\bibitem{BCDGZ} Bracci, F., Contreras, M.D., D\'{\i}az-Madrigal, S., Gaussier, H., Zimmer, A.:
Asymptotic behavior of orbits of holomorphic semigroups. J. Math. Pures Appl. (9) {\bf133} (2020), 263--286.

\bibitem{BCD-Book} Bracci, F., Contreras, M.D., D\'{\i}az-Madrigal, S.: \textsl{Continuous semigroups of holomorphic functions in the unit disc.} Springer Monographs in Mathematics, 2020.

\bibitem{ConDia05a} Contreras, M.D., D\'{\i}az-Madrigal, S.: Analytic flows on the unit disk: angular derivatives and boundary fixed points. Pacific J. Math. \textbf{222}  (2005), 253--286.

\bibitem{CDG-Slope} Contreras, M.D., D\'{\i}az-Madrigal, S., Gumenyuk, P.: Slope problem for trajectories of holomorphic semigroups in the unit disk. Comput. Methods Funct. Theory {\bf 15} (2015), no.~1, 117--124. 

\bibitem{Cow81} Cowen, C.C.:  Iteration and the solution of functional equations for functions analytic in the unit disk. Trans. Amer. Math. Soc. \textbf{265} (1981), 69--95.

\bibitem{EKRS} Elin M., Khavinson D., Reich S., and Shoikhet D.: Linearization models for parabolic dynamical systems via Abel’s functional equation. Ann. Acad. Sci. Fen. Math. {\bf 35} (2010), 439–472.

\bibitem{Goluzin} Goluzin, G.M.: \textsl{Geometric theory of functions of a complex variable}. Amer. Math. Soc., Providence, R.I., 1969. MR0247039 (Translated from G. M. Goluzin, {\it Geometrical theory of functions of a complex variable} (Russian), Second edition, Izdat. ``Nauka'', Moscow, 1966.)

\bibitem{Kelgiannis} Kelgiannis, G.: Trajectories of semigroups of holomorphic functions and harmonic measure. J. Math. Anal. Appl. {\bf 474} (2019), no.~2, 1364--1374. 

\bibitem{Pombook75} Pommerenke, Ch.: \textsl{Univalent functions. With a chapter on quadratic differentials by Gerd Jensen}. Vandenhoeck \& Ruprecht, G\"{o}ttingen, 1975. 

\bibitem{Pombook92} Pommerenke, Ch.: \textsl{Boundary behaviour of conformal maps}, Grundlehren der Mathematischen Wissenschaften, 299, Springer-Verlag, Berlin, 1992. 

\bibitem{Shobook01} Shoikhet, D.: \textsl{Semigroups in geometrical function theory}. Kluwer Academic Publishers, Dordrecht, 2001.

\bibitem{Sis} Siskakis, A.G.: \textsl{Semigroups of composition operators and the Ces\`aro operator on
$H^p(D)$}. Ph.D. thesis, University of Illinois, 1985.

\end{thebibliography}
\end{document}